# STATIONARY DISTRIBUTIONS OF A MODEL OF SYMPATRIC SPECIATION[1]

By Feng Yu

*University of Oxford*

This paper deals with a model of sympatric speciation, that is, speciation in the absence of geographical separation, originally proposed by U. Dieckmann and M. Doebeli in 1999. We modify their original model to obtain a Fleming–Viot type model and study its stationary distribution. We show that speciation may occur, that is, the stationary distribution puts most of the mass on a configuration that does not concentrate on the phenotype with maximum carrying capacity, if competition between phenotypes is intense enough. Conversely, if competition between phenotypes is not intense, then speciation will not occur and most of the population will have the phenotype with the highest carrying capacity. The length of time it takes speciation to occur also has a delicate dependence on the mutation parameter, and the exact shape of the carrying capacity function and the competition kernel.

**1. Introduction.** Understanding speciation is one of the great problems in the field of evolution. According to Mayr [9], speciation means the splitting of a single species into several, that is, the multiplication of species. It is believed that many species originated through geographically isolated populations of the same ancestral species. This phenomenon is relatively easy to understand. In contrast, sympatric speciation, in which new species arise without geographical isolation, is theoretically much more difficult. In this work, we take the recent work in Dieckmann and Doebeli [4] on sympatric speciation as a basis, and try to develop a model that captures the most important aspects of their model and yet is also amenable to rigorous mathematical analysis. In Section 1.1, we briefly describe the Dieckmann–Doebeli

Received August 2005; revised October 2006.

[1]Supported by the Department of Mathematics and a Killam Fellowship while at the University of British Columbia, Vancouver, Canada, and EPSRC Grant GR/T19537 while at the University of Oxford.

*AMS 2000 subject classifications.* Primary 92D15; secondary 60J25, 60J27, 92D25.

*Key words and phrases.* Fleming–Viot processes, sympatric speciation, population dynamics, competition model, evolutionary biology.

model of sympatric speciation. Their original model is very difficult to study, so in Section 1.2, we present a simplified model that retains almost exactly the fitness function found in the original Dieckmann–Doebeli model, and perform some nonrigorous analysis that illustrates the delicate dependence of transitory behavior on the exact form of the fitness function. In Section 2, we present our main model, a Fleming–Viot model with strong selection and a fitness function that retains the key features of the original Dieckmann–Doebeli model. The advantage of using a Fleming–Viot model is that one can write down the stationary distribution quite explicitly, and stationary or long-term behavior is usually easier to study than transitory ones. It turns out that the stationary distribution concentrates more and more mass near its global maximum as the population size becomes larger, thus identifying the global maximum gives a strong indication of the kind of configuration eventually taken up by the population. The main results are given toward the end of Section 2, along with some discussion of these results. The rest of the paper, Section 3, is devoted to proofs of various results on local and global maxima of the stationary distribution of the Fleming–Viot model introduced in Section 2.

1.1. *The Dieckmann–Doebeli model.* Dieckmann and Doebeli [4] proposed a general model for sympatric speciation, for both asexual and sexual populations. We will briefly describe their model for the asexual population, since this is the model we study in this work. Their sexual model is naturally more complicated than the asexual model, but the two models have similar behavior. In their asexual model, each individual in the population is assumed to have a quantitative character (phenotype) $x \in \mathbb{R}$ determining how effectively this individual can make use of resources in the surrounding environment. A typical example is the beak size of a certain bird species, which determines the size of seeds that can be consumed by an individual bird. The function $K : \mathbb{R} \to \mathbb{R}^+$ (carrying capacity) is associated with the surrounding environment, where $K_x$ denotes the number of individuals of phenotype $x$ that can be supported by the environment. For example, since birds with small beak size (say $x_1$) are more adapted to eating small seeds than birds with large beak size (say $x_2$, $x_2 > x_1$), $K_{x_1}$ will be larger than $K_{x_2}$ if the surrounding environment produces more small seeds than large seeds. In the Dieckmann–Doebeli model, $K_x$ is taken to be $c\exp(-(x-\hat{x})^2/(2\sigma_K^2))$. Moreover, every pair of individuals compete at an intensity determined by the phenotypical distance of these two individuals. More specifically, an individual of phenotype $x_1$ competes with an individual of phenotype $x_2$ at intensity $C_{x_1-x_2}$, where $C_x = \exp(-x^2/(2\sigma_C^2))$. Therefore each individual in the population interacts with the environment via the carrying capacity $K$, and interacts with the population via the competition kernel $C$.

Let $N_x(t)$ denote the number of individuals with phenotype $x$ at time $t$. At any time, an individual of phenotype $x$ gives birth at a constant rate, and dies at a rate proportional to $(C * N.(t))_x / K_x$, that is, inversely proportional to the $x$-carrying capacity, but proportional to the intensity of competition exerted by the population on phenotype $x$, the numerator $(C * N.(t))_x = \int C_{x-y} N_y(t)\, dy$ being how much competition (from every individual in the population) individuals with phenotype $x$ suffer. In addition, every time an individual gives birth, there is a small probability that a mutation occurs and the phenotype of the offspring is different from that of the parent; in this case, the phenotypical distance between the offspring and the parent is then random and assumed to have a Gaussian distribution.

Since the number of individuals of a certain phenotype increases via the birth mechanism at a linear rate, but decreases via the death mechanism at a quadratic rate, extinction of all phenotypes will occur in finite time with probability one, that is, $N \equiv 0$ eventually. For large initial populations, however, extinction will happen far enough into the future that interesting behavior does arise before the population becomes extinct.

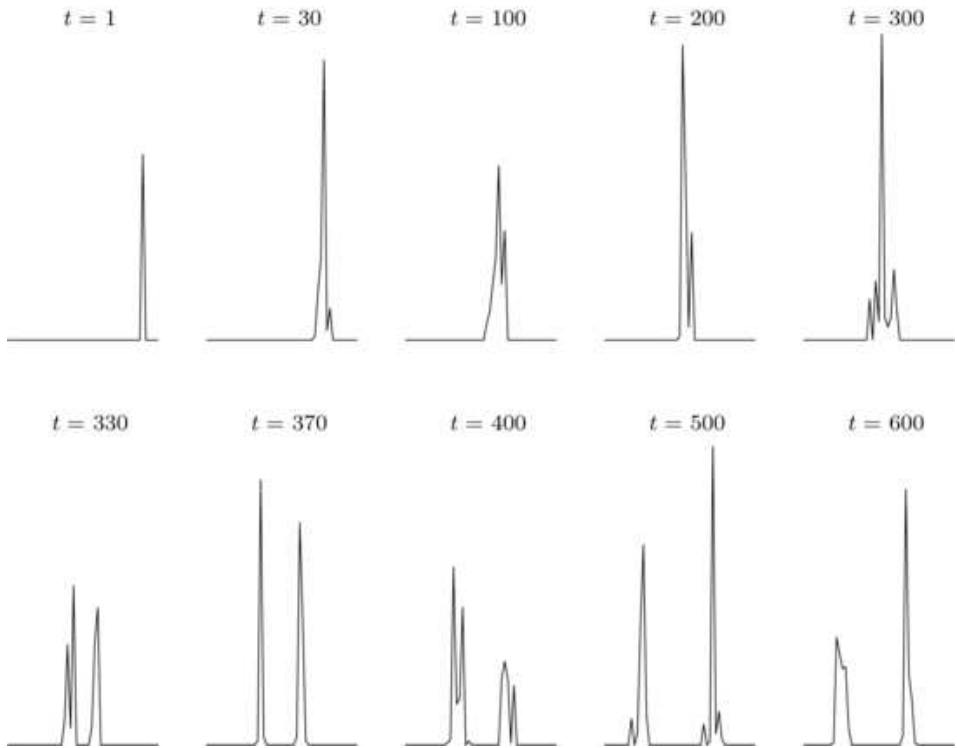

FIG. 1. *Simulation of the Dieckmann–Doebeli model with $E = [-50, 50] \cap \mathbb{Z}$, $\sigma_K = \sqrt{1000}$, $\sigma_C = \sqrt{600}$ and mutation happening to 1.5% of the births.*

Monte–Carlo simulations, shown in Figure 1, give a fairly good idea of the behavior of the Dieckmann–Doebeli model for asexual populations. If the initial population is monomorphic ($t = 1$ in Figure 1), that is, concentrated near a certain phenotype $x_0$ ($N.(0)/\sum_x N_x(0) \approx \delta_{x_0}$), then the entire population first moves ($t = 30, 100, 200$ in Figure 1) toward $\hat{x}$, the phenotype with maximum carrying capacity. If $\sigma_C > \sigma_K$ (this includes the case $\sigma_C = \infty$, i.e., equal competition between all phenotypes), then the population stabilizes near phenotype $\hat{x}$. But if $\sigma_C < \sigma_K$, then the monomorphic population concentrated at phenotype $\hat{x}$ splits into two groups, one group concentrating on a phenotype $< \hat{x}$, while the other concentrating on a phenotype $> \hat{x}$ ($t = 330, 370, 400, 500$ in Figure 1). In the latter case, one can say that one species has evolved into two distinct species. The variance of the Gaussian distribution used in the mutation kernel affects how different phenotypically the offspring can be from the parent, and seems to affect the speed of evolution, but not the configuration eventually taken up by the population.

1.2. *A conditioned Dieckmann–Doebeli model.* As noted in the very first paragraph, the Dieckmann–Doebeli model for asexual populations is very difficult to study. One reason for this difficulty is because the number of individuals can fluctuate with time. As mentioned before, since the birth rate is linear but the death rate is quadratic, extinction will occur in finite time with probability one, which makes it somewhat meaningless to analyze the stationary or long-term behavior of the system. The modification we apply to the Dieckmann–Doebeli model is to assume constant population size $N$, reflecting a constant carrying capacity of the overall population, and define a Wright–Fisher type model (for a definition of Wright–Fisher model and its relationship with Fleming–Viot models, see [6]) with fitness functions chosen to retain key ingredients of the original Dieckmann–Doebeli model. In contrast to the continuous-time nature of the original Dieckmann–Doebeli model, the modified model is discrete time. Because the number of individuals remains constant, analyzing the behavior of the population is equivalent to analyzing the empirical distribution

$$\pi^N = \frac{1}{N} \sum_{n=1}^{N} \delta_{x_n},$$

where $x_n$, $n = 1, \ldots, N$, denotes the phenotype of the $n$th individual in a population of size $N$ and $\delta_x$ is the measure that puts unit mass at phenotype $x$.

Before we describe our choice of fitness functions, we briefly describe the concepts of *fitness* and *selection*. Selection occurs when individuals of different genotypes leave different numbers of offspring because their probabilities of surviving to reproductive age are different (see [1]). If we define fitness to be a measure of how likely a particular individual produces offspring that will survive to reproductive age, then individuals with higher fitness should have higher probability of being selected for reproduction. Along these lines, it is natural to define fitness of a phenotype as the difference between the birth rate and the death rate of individuals of this phenotype, therefore it is also natural to require the fitness function to be bounded.

The key feature of the Dieckmann–Doebeli model is that each individual has a fitness that depends on both the carrying capacity associated with its phenotype and the configuration of the entire population. More specifically, the fitness of a phenotype $x$ is an increasing function of $K_x$, the carrying capacity, but a decreasing function of $(C * N)_x$, the competition it suffers. Here $N_x$ is the number of individuals of phenotype $x$. With this in mind, we propose the following two fitness functions:

$$W_x^{(1)}(\pi) = 0 \vee \left(1 - \frac{\sum_z C_{x-z}\pi_z}{K_x}\right),$$

$$W_x^{(2)}(\pi) = \frac{K_x}{\sum_z C_{x-z}\pi_z}.$$

Each of the two fitness function defined above is an increasing function of $K_x$ and a decreasing function of $(C * \pi)_x$. $W^{(1)}$ resembles more closely the original Dieckmann–Doebeli model, but it has the disadvantage of being in a more complicated form than $W^{(2)}$ and it is also not differentiable. Our simplified discrete-time and discrete-space Dieckmann–Doebeli model is as follows:

- At every time step $t \in \mathbb{Z}^+$, the entire population is replaced by a new population of $N$ individuals, each chosen independently according to the distribution $p.(t, \pi^N)$:

$$p_x(t, \pi^N) = \sum_y A(y, x) \frac{\pi_y^N(t) W_y(\pi^N(t))}{\sum_z \pi_z^N(t) W_z(\pi^N(t))}$$

where the denominator $\sum_z \pi_z^N(t) W_z(\pi^N(t))$ is simply the normalization factor such that $\sum_x p_x(t, \pi^N) = 1$ and

1. $E = [-L, L] \cap \mathbb{Z}$ is the phenotype space, and $\pi^N \in \mathcal{P}(E)$ is a probability measure on $E$,
2. $K : E \to [0, 1]$ is the carrying capacity, and $C : \mathbb{Z} \to \mathbb{R}^+$ is the competition kernel,
3. $W_x(\pi)$ is the fitness of phenotype $x$ in a population with empirical distribution $\pi$ (sometimes we notationally suppress the dependence on $\pi$), and $W = W^{(1)}$ or $W = W^{(2)}$,

4. $A$ is a Markov transition matrix associated with mutation, with $A(y,x)$ denoting the probability of an individual of phenotype $y$ mutating to an individual of phenotype $x$.

By Theorem 1 in [3], $\{\pi_t^N, t \in [0,T]\} \Rightarrow \{\pi_t, t \in [0,T]\}$ as $N \to \infty$, where $\Rightarrow$ denotes weak convergence and $\pi_t$ evolves according to the following *deterministic* dynamical system:

$$\pi_x(t+1) = \sum_y A(y,x) \frac{\pi_y(t) W_y(\pi(t))}{\sum_z \pi_z(t) W_z(\pi(t))}. \tag{1}$$

Analyzing the dynamical system (1) is still not easy, partly because it is of a complicated form that is nonlinear in $\pi$, and we cannot find any Lyapunov function [8] that associates with (1). Simulations of (1), however, seem to display some interesting behavior, which we will describe after carrying out some nonrigorous analysis of (1).

Without mutation, any phenotype $x$ with $\pi_x = 0$ at any time $\tau$ will stay 0 for all $t \geq \tau$. Mutation enables individuals of phenotype $x$ to be born in future generations even if there are no individuals of phenotype $x$ in the present generation. But if we start with a polymorphic initial measure, that is, $\pi_x(0) \neq 0$ for all $x$, then adding small mutation to the system should not cause significant changes in the behavior of (1). Therefore we assume that $A = I$ and $\pi(0)$ is polymorphic. In this case, (1) can be simplified to

$$\pi_x(t+1) = \frac{\pi_x(t) W_x(\pi(t))}{\sum_z \pi_z(t) W_z(\pi(t))}.$$

Thus if $A = I$, then $\hat{\pi}$ is a stationary distribution of (1) if and only if

$$\hat{\pi}_x = \frac{1}{c} \hat{\pi}_x W_x(\hat{\pi}) \tag{2}$$

for some constant $c$. Condition (2) is equivalent to

$$W_x(\hat{\pi}) = c \quad \text{for all } x \text{ where } \hat{\pi}(x) \neq 0. \tag{3}$$

Let $K$ and $C$ be in the form considered by Dieckmann and Doebeli, that is, $K_x = \exp(-x^2/2\sigma_K^2)$ and $C_x = \exp(-x^2/2\sigma_C^2)$. If $W = W^{(2)}$, then condition (3) means that

$$K_x = c(C * \hat{\pi})(x) \quad \text{for all } x \text{ where } \hat{\pi}(x) \neq 0,$$

which seems to indicate that if $\sigma_C < \sigma_K$, then $\hat{\pi}$ should be close to $\mathcal{N}(0, \sigma_K^2 - \sigma_C^2)$. On the other hand, if $W = W^{(1)}$, then $\hat{\pi}$ is a stationary distribution if $1 - (\sum_z C_{x-z} \hat{\pi}_z)/K_x$ is a strictly positive constant. Notice that if $K$ and $C$ are both Gaussian-shaped with $K_0 = C_0 = 1$ then $\hat{\pi} = \mathcal{N}(0, \sigma_K^2 - \sigma_C^2)$ makes $1 - (\sum_z C_{x-z} \hat{\pi}_z)/K_x$ constant; furthermore, this constant is strictly positive since $(C * \hat{\pi})(0) < K_0 = 1$ if $\sigma_C < \sigma_K$.

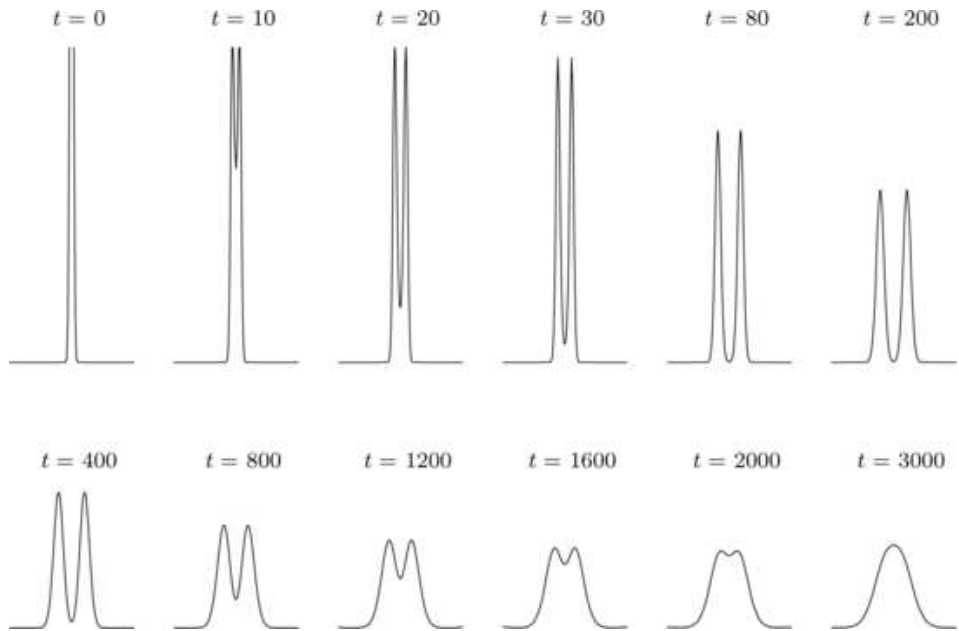

FIG. 2. *Simulation of* (1) *with* $E = [-149, 149] \cap \mathbb{Z}$, $\sigma_K = 60$, $\sigma_C = 55$ *and* $W = W^{(1)}$.

Therefore for both $W^{(1)}$ and $W^{(2)}$, assuming Gaussian competition and carrying capacity kernels, the dynamical system (1) should have Gaussian-shaped stationary distributions if $\sigma_C < \sigma_K$. In simulations carried out by Dieckmann and Doebeli [4], however, $\sigma_C < \sigma_K$ is the case that leads to speciation, that is, the stationary distribution supposedly has two sharp well-separated peaks, which contradicts the analysis carried out in the previous paragraph. Simulations of (1) with $W = W^{(1)}$, shown in Figure 2, reveal that if $\pi(0) \approx \delta_0$, initially the population does split into two groups and begins to move apart, but as $t \to \infty$, the empirical measure converges to a Gaussian-shaped hump. This suggests the possibility that in the original Dieckmann–Doebeli model, conditioning on the population surviving long enough for convergence to stationarity to occur (recall that in the original Dieckmann–Doebeli model, extinction occurs in finite time), speciation is also a transitory phenomenon, rather than a stationary phenomenon. Simulations of (1) with $W = W^{(2)}$, shown in Figure 3, do not even display transitory speciation behavior. Instead, the initial spike at 0 simply widens to a Gaussian hump centered at 0. Hence the particular form of the dependence on $K$ and $C * \pi$ seems to affect whether or not speciation occurs.

From the simulations and nonrigorous analysis above, it seems that the dynamical system in (1) does not have a bimodal stationary distribution if both $K$ and $C$ are taken to be Gaussian-shaped. But we would also like to point out that if the shape of $K$ or $C$ were changed just a bit, for example, making $K$ a bit more "flat" by taking $K_x = \exp(-x^{2-\varepsilon}/2\sigma_K^2)$, simulations then display stationary distributions that have two (or even more) modes. Thus the shape of the stationary distributions of the Dieckmann–Doebeli model seems to have a very delicate dependence on the choice of $K$ and $C$. We also speculate that since the original Dieckmann–Doebeli model has a fluctuating population size whereas our nonrigorous analysis only applies to a model with fixed population size, this small difference may also disturb the long-time behavior of the model enough that a Gaussian hump does not appear with high probability before the population becomes extinct.

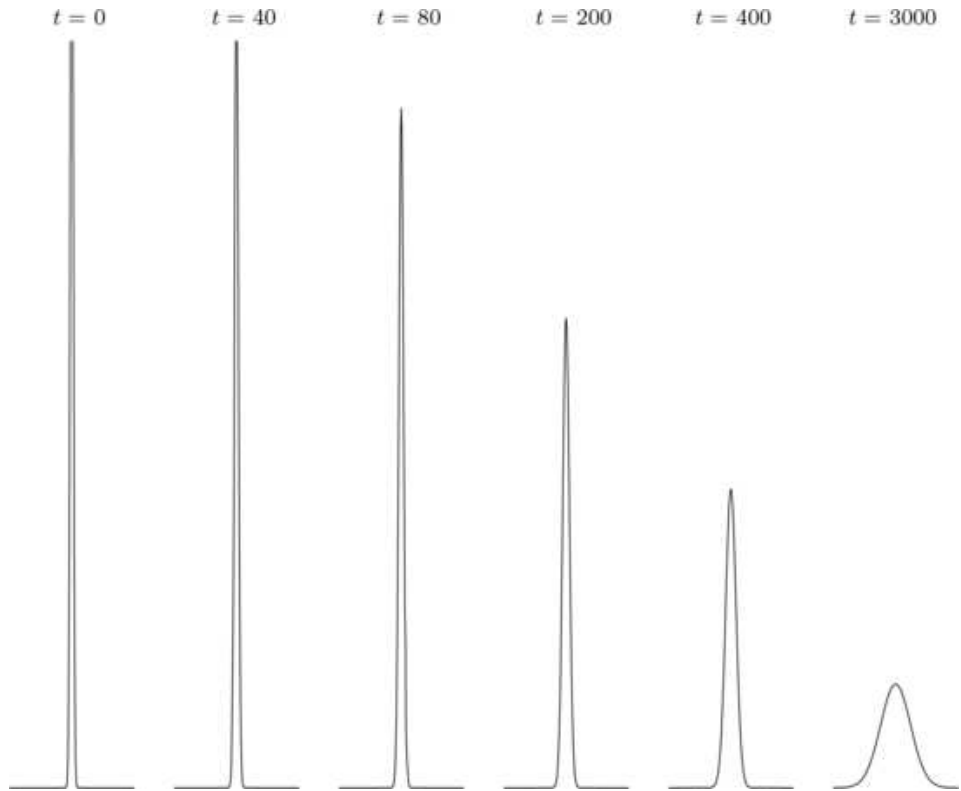

Fig. 3. *Simulation of* (1) *with* $E = [-149, 149] \cap \mathbb{Z}$, $\sigma_K = 60$, $\sigma_C = 55$ *and* $W = W^{(2)}$.

In case of this transitory behavior, we may still say speciation has occurred. A constant carrying capacity function is only an approximation of what actually happens in nature, where the environment a species lives in can change quite drastically over a long period of time. By assuming a carrying capacity function that does not change over time, we are essentially studying what can happen to a single species over time lengths during which this approximation is reasonable.

We also have a rigorous result for this discrete-time model, for the simplest case of $K$ and $C$ both rectangular, that is, $K_x = 1_{\{|x| \leq L\}}$ and $C_x = 1_{\{|x| \leq M\}}$ for some integers $L$ and $M$. Theorem A.0.14 from [10] says that if $\nu^n$ is a convergent sequence of symmetric stationary distributions for the conditioned Dieckmann–Doebeli model where the mutation matrix $A$ corresponds with a convolution kernel $\mu^n \delta_{-1} + (1 - 2\mu^n)\delta_0 + \mu^n \delta_1$ then $\nu^n([-(M-L+1), M-L+1]) \to 0$ as $\mu^n \to 0$; in words, the mass in the middle gets very small as the mutation parameter approaches zero, hence there exist bimodal stationary distributions if $\mu^n$ is sufficiently small.

In the next section, we introduce the Fleming–Viot model that we study for the rest of the work. It is a continuous-time model that approximates a Moran particle system. The main advantage of a Fleming–Viot type model is that if the fitness function is chosen to be a quadratic form in $\pi$, then the exact form of the stationary distribution is known in the literature [7].

**2. The Fleming–Viot model and main results.** We work on the phenotype space $E = [-L, L] \cap \mathbb{Z}$. Sometimes we refer to a phenotype as a site in $E$. Let

$$\Delta = \left\{ (\pi_{-L}, \ldots, \pi_0, \ldots, \pi_L) : \pi_i \geq 0 \ \forall i \text{ and } \sum_{i=-L}^{L} \pi_i = 1 \right\}$$

be the space of probability measures on $E$, that is, $\Delta = \mathcal{P}(E)$. Members of $\Delta$ are usually denoted by $\pi$, $\hat{\pi}$, $\pi^N$, and so on. We endow $\Delta$ with the following metric:

$$d(\hat{\pi}, \tilde{\pi}) = \max_x |\hat{\pi}(x) - \tilde{\pi}(x)|.$$

We assume a monomorphic initial condition, that is, $\pi_0 = \delta_x$ for some $x \in E$ (in fact, we take $x = 0$ mostly).

Recalling that the essential ingredient of the original Dieckmann–Doebeli model is that the fitness function is an increasing function of $K_x$ and a decreasing function of $(C * \pi)_x$, we define fitness $m_x(\pi)$ and mean fitness $\overline{m}_\pi$ for our Fleming–Viot model (for a precise definition the Fleming–Viot process, see [5] or [2]) to have the following form:

$$m_x(\pi) = K_x \sum_z B_{x-z} K_z \pi_z,$$

(4)

$$\overline{m}_\pi = \sum_x \pi_x m_x,$$

where the "cooperation" kernel $B$ can be taken to be $1 - C$. We assume $B$ is symmetric. In the original Dieckmann–Doebeli model, pairs of individuals with small phenotypical distance *compete* at a higher intensity than pairs of individuals with large phenotypical distance; in our model, pairs of individuals with small phenotypical distance *cooperate* at a lower intensity than pairs of individuals with large phenotypical distance.

The term $B_{x-z}$ in the definition of $m_x(\pi)$ above can be thought of as a measure of how *inefficiently* an individual of phenotype $z$ makes use of resources of type $x$, that is, the type that best suit individuals of phenotype $x$. For example, if individuals of phenotype $z$ cannot makes use of resources of type $x$ at all, that is, $B_{x-z} = 1$, then they contribute to an increase to the fitness of individuals of phenotype $x$, since this type $z$ individual will not compete with type $x$ individuals. On the other hand, if individuals of phenotype $z$ makes perfect use of resources of type $x$, that is, $B_{x-z} = 0$, then these individuals contribute no increase to the fitness of individuals of phenotype $x$. From the point of view of a particular individual of phenotype $x$, he "prefers" (if he is selfish) all other individuals in the population to be of phenotypes $z$ with $B_{x-z} = 1$, so that no other individual can make use of resources for which he is best adapted. Thus the term "cooperation" is somewhat misleading, since individuals with different phenotypes do not really cooperate with each other. Nevertheless, we use "cooperation" and "competition" to describe the effect of individuals of a certain phenotype on individuals of another phenotype out of convenience. If $B_z = 0$, then we say phenotypes separated by distance $z$ do not cooperate at all (i.e., compete at full intensity), and if $B_z = 1$, we say they cooperate at full intensity (i.e., do not compete at all).

2.1. *The model.* Let $K : E \to [0,1]$ be the carrying capacity function, and $B : \mathbb{Z} \to [0,1]$ be the cooperation kernel. We assume $B$ is symmetric. We define

$$\begin{aligned}(5)\quad \mathcal{G} = &\sum_{x=-L}^{L} \left[\frac{\mu}{2}(1 - (2L+1)\pi_x) + \sigma\pi_x(m_x(\pi) - \overline{m}_\pi)\right]\frac{\partial}{\partial \pi_x} \\ &+ \frac{1}{N}\sum_{x,y=-L}^{L} \pi_x(\delta_{xy} - \pi_y)\frac{\partial^2}{\partial \pi_x \partial \pi_y}\end{aligned}$$

to be the generator of our Fleming–Viot process with selection and mutation, where $\delta_{xy} = 1$ if $x = y$ and $= 0$ otherwise, and the fitness of site $x$ in a population with distribution $\pi$ and the mean fitness of the population $\overline{m}_\pi$ are defined in (4). A Fleming–Viot process with finitely many types is also known as a Wright–Fisher diffusion (see [2]), but to stress the continuous time nature and avoid confusion with the discrete time Wright–Fisher model, we still refer to our model as a Fleming–Viot process, which is a special case of the Fleming–Viot process with selection as described in Chapter 10.1.1 of [2].

In (5), the terms that correspond with the effect of selection and replacement sampling are the following:

$$\mathcal{G}^S = \sigma \sum_{x=-L}^{L} \pi_x(m_x(\pi) - \overline{m}_\pi)\frac{\partial}{\partial \pi_x} + \frac{1}{N} \sum_{x,y=-L}^{L} \pi_x(\delta_{xy} - \pi_y)\frac{\partial^2}{\partial \pi_x \partial \pi_y}.$$

$\mathcal{G}^S$ approximates the following Moran particle system (see page 26 of [2] for a precise definition of Moran particle systems) with a population of $N$ individuals undergoing *strong* selection for suitably small $\sigma$ (e.g., $\sigma \leq 1/2$ if $K \leq 1$):

- $\pi_x^N$ decreases by $1/N$ and $\pi_y^N$ increases by $1/N$ at rate $\frac{N}{2}\pi_x^N(\frac{1}{2}+\sigma(m_y(\pi^N) - m_x(\pi^N)))\pi_y^N$.

To see this, we expand the generator $\mathcal{G}^{S,N}$ for the particle system above for smooth and compactly supported $f(\pi_{-L}, \ldots, \pi_L): \mathbb{R}^{2L+1} \to \mathbb{R}$:

$$\mathcal{G}^{S,N} f(\pi^N) = \sum_{x,y=-L}^{L} \left[ f\left(\pi^N - \frac{1}{N}\delta_x + \frac{1}{N}\delta_y\right) - f(\pi^N) \right] N\pi_x^N$$

$$\times \left(\frac{1}{4} + \frac{\sigma}{2}(m_y(\pi^N) - m_x(\pi^N))\right)\pi_y^N$$

$$= \sum_{x,y=-L}^{L} \left[\frac{\partial f(\pi^N)}{\partial \pi_y} - \frac{\partial f(\pi^N)}{\partial \pi_x} + \frac{1}{2N}\frac{\partial^2 f(\pi^N)}{\partial \pi_y^2} + \frac{1}{2N}\frac{\partial^2 f(\pi^N)}{\partial \pi_x^2}\right.$$

$$\left.- \frac{1}{N}\frac{\partial^2 f(\pi^N)}{\partial \pi_x \partial \pi_y} + O(1/N^2)\right]$$

$$\times \pi_x^N\left(1 + \frac{\sigma}{2}(m_y(\pi^N) - m_x(\pi^N))\right)\pi_y^N$$

$$= \sigma \sum_{x=-L}^{L} \pi_x^N(m_x(\pi^N) - \overline{m}_{\pi^N})\frac{\partial f(\pi^N)}{\partial \pi_x}$$

$$+ \frac{1}{N} \sum_{x,y=-L}^{L} \pi_x^N(\delta_{xy} - \pi_y^N)\frac{\partial^2 f(\pi^N)}{\partial \pi_x \partial \pi_y} + O(1/N^2).$$

Therefore the generators of the particle system and the Fleming–Viot process without mutation, $\mathcal{G}^{S,N}$ and $\mathcal{G}^S$ respectively, agree up to $O(1/N)$, and as $N \to \infty$, the stochastic process associated with them both converge to the solution of the following system of deterministic ordinary differential equation (ODEs).

(6) $$\partial_t \pi_x = \sigma \pi_x(m_x(\pi) - \overline{m}_\pi).$$

In the mutation component of $\mathcal{G}$, $\sum_{x=-L}^{L} \mu(1-(2L+1)\pi_x)\partial/\partial \pi_x$, we use the simplifying assumption that the mutation is symmetric, that is, the rate $\mu_{xy} = \mu_y$ at which phenotype $x$ mutates to phenotype $y$ depends on $y$ only, and furthermore $\mu_y = \mu$ is constant in $y$; the latter assumption makes the proofs a bit cleaner. In the original Dieckmann–Doebeli model, the variance of the mutation kernel only affects the speed of evolution, not the eventual configuration taken up by the population. Therefore the assumption of symmetric mutation should not affect the stationary behavior of the process a great deal, and it is precisely this assumption that enables one to write down the unique stationary distribution for the Fleming–Viot process, as well as a Lyapunov function for its infinite population limit. Furthermore, the mutation component of $\mathcal{G}$ is an approximation of the $N$-particle system that undergoes the following:

- $\pi_x^N$ decreases by $1/N$ and $\pi_y^N$ increases by $1/N$ at rate $\frac{N}{2}\mu\pi_x^N$, with an error term of $O(\mu/N)$. We can expand the generator $\mathcal{G}^{M,N}$ of the particle system above:

$$\mathcal{G}^{M,N} f(\pi^N) = \sum_{x,y=-L}^{L} \left[ f\left(\pi^N - \frac{1}{N}\delta_x + \frac{1}{N}\delta_y\right) - f(\pi^N) \right] \frac{N}{2}\mu\pi_x^N$$

$$= \frac{\mu}{2} \sum_{x=-L}^{L} (1 - (2L+1)\pi_x^N) \frac{\partial f(\pi^N)}{\partial \pi_x^N}$$

$$+ \frac{\mu}{4N} \sum_{x,y=-L}^{L} \pi_x^N \left[ \frac{\partial^2 f(\pi^N)}{\partial \pi_y^{N2}} + \frac{\partial^2 f(\pi^N)}{\partial \pi_x^{N2}} - 2\frac{\partial^2 f(\pi^N)}{\partial \pi_x^N \partial \pi_y^N} \right]$$

$$+ O(1/N^2),$$

which has a rather messy noise term (the term involving second derivatives of $f$). The interesting cases are those with small $\mu$, and we only retain the drift term (terms involving first derivatives) in the expansion of $\mathcal{G}^{M,N}$. Combining this with the drift and noise terms in the expansion of $\mathcal{G}^{S,N}$, we obtain the generator $\mathcal{G}$. As $N \to \infty$, the process with generator $\mathcal{G}$ converges to the solution of the following system of deterministic ODEs:

$$\partial_t \pi_x = \sigma \pi_x (m_x(\pi) - \overline{m}_\pi) + \frac{\mu}{2}(1 - (2L+1)\pi_x). \tag{7}$$

One can apply Theorem 1.6.1 from [5] to establish this convergence.

The generator $\mathcal{G}$ is of the form defined in Lemma 4.1 from [7] if one speeds up time by $N/2$ in (5), and a direct application of that result implies the following result:

PROPOSITION 2.1. *For the Fleming–Viot process with generator $\mathcal{G}$,*

$$\nu^N(d\pi) = C \left( \prod_{x=-L}^{L} \pi_x \right)^{(N/2)\mu - 1} e^{(N/2)\overline{m}_\pi} d\pi_{-L} \cdots d\pi_L$$

*is the unique stationary distribution, where $C$ is the normalizing constant such that $\nu^N$ is a probability measure on $\Delta$.*

We define $\tilde{\mu} = \mu - \frac{2}{N}$, which we assume to be positive, and write

$$\nu^N(d\pi) = C \exp\left\{ \frac{N}{2} \left( \overline{m}_\pi + \left(\mu - \frac{2}{N}\right) \sum_{x=-L}^{L} \log \pi_x \right) \right\} d\pi_{-L} \cdots d\pi_L.$$

As $N \to \infty$, we expect $\nu^N$ to concentrate more and more on the configuration that maximizes

$$(8) \qquad V_\pi = \overline{m}_\pi + \tilde{\mu} \sum_{x=-L}^{L} \log \pi_x.$$

One can imagine a scenario where initially all birds in the population have beaks that specialize in eating seeds of say size 5, 5 being the most common size in the forest. As time passes, the selection part $\mathcal{G}$, $\sum_{x=-L}^{L} \sigma \pi_x (m_x(\pi) - \overline{m}_\pi) \partial/\partial \pi_x$, moves the population toward a fitter configuration, since the mean fitness $\overline{m}_\pi$ is a Lyapunov function of the dynamical system (6). [A proof of a slightly more general statement can be found in Lemma 2.4(a).] If the forest produces nearly as many seeds of size 4 and 6 as seeds of size 5, but the birds can really just eat one size of seeds (e.g., if a bird's beak specializes in seeds of size 5, then it is very bad at eating seeds of size 4 or 6), then it is quite possible that the population as a whole does better, that is, is more fit on average, if half the birds specialize in seeds of size 4, while the other half in seeds of size 6. This way, even though each bird has to spend slightly more effort to find seeds that suits her (since $K_5$ is slightly larger than $K_4$ or $K_6$), she only competes with half the population.

Proposition 2.2 below says that as the population size becomes large and the mutation parameter becomes small, the stationary distribution $\nu^N$ focuses more and more on the configuration that achieves the maximum fitness. This configuration is also the one that the population spends the most time in. Going back to the example in the last paragraph, if it can be verified that the configuration where half the birds specialize in seeds of size 4 while the other half in seeds of size 6 maximizes fitness, then starting from an initial population where all birds specialize in seeds of size 5, the population will eventually drift to the fitter bimodal configuration. In this case, we can say that speciation has occurred. Because of the stochastic nature of the model, eventually the population will leave this maximally fit configuration and enter some less fit configuration. But this may not happen for a long time, after which the validity of the approximation that carrying capacity is constant over time may no longer be valid.

We broadly say that speciation is likely to occur eventually if the population configuration that achieves maximum $V_\pi$ has *significant* mass at phenotype(s) different from the original one (the original phenotype may or may not die out as a result of speciation). In general, however, it is difficult to identify the configuration that (globally) achieves maximum $V_\pi$, or even to verify that a certain configuration achieves it, due to the non-concave nature of the $V_\pi$ and even $\overline{m}_\pi$. Studying local maxima of $V_\pi$ then becomes useful, where one can exclude certain classes of configuration from candidates for the global maximum of $V_\pi$. For example, if one can exclude configuration close to $\delta_0$ as a local maximum, then the global maximum will have significant mass at sites other than 0, and we may also say speciation will eventually occur in this case.

Mutation effects alone produces individuals of all phenotypes in $E$, so the word "significant" in our definition of speciation is taken to mean a number that does not go to 0 as the mutation parameter $\mu$ goes to 0. Since mutation just distributes mass evenly to all sites in $E$, a large $\mu$ obscures the effects of selection, thus we are mainly interested in small $\mu$. We attempt to bound mass at various phenotypes away from 0 as $\mu \to 0$. The case of $\mu = 0$ is not interesting, since with initial condition $\delta_0$ the configuration will then remain monomorphic and no speciation can ever occur.

Propositions 2.2 and 2.3 below relate local/global maximum of $\overline{m}_\pi$ to those of $V_\pi$ when $\mu$ is small. In particular, as can be expected, they are quite close to each other when $\mu$ is small.

PROPOSITION 2.2. *Let $\{\tilde{\pi}_1, \tilde{\pi}_2, \ldots, \tilde{\pi}_k\}$ be the finite set that consists of all global maxima of $\overline{m}_\pi$ for $\pi \in \Delta$. For any sufficiently small $\varepsilon > 0$, we can pick $\mu$ small and $N$ large, such that*

$$\nu^N\left(\bigcup_{i=1}^k \mathrm{Ball}(\tilde{\pi}_i, \varepsilon)\right) > 1 - \varepsilon,$$

*where $\mathrm{Ball}(\pi, \varepsilon)$ denotes the intersection of $\Delta$ and the ball of radius $\varepsilon$ centered at $\pi$.*

PROPOSITION 2.3. *Let $\tilde{\pi}$ be the* unique *local maximum of $\overline{m}_\pi$ for $\pi \in \Delta$ in a small neighborhood of $\tilde{\pi}$. If $\tilde{\mu}$ is sufficiently small, then $V_\pi$ as defined by (8) has a local maximum in a small neighborhood of $\tilde{\pi}$.*

LEMMA 2.4. *We have (a) $V_\pi : \overset{\circ}{\Delta} \to \mathbb{R}$ is a Lyapunov function for the dynamical system*

$$\partial_t \pi_x = \pi_x\left(m_x - \overline{m}_\pi + \frac{\tilde{\mu}}{2}\left(\frac{1}{\pi_x} - (2L+1)\right)\right), \tag{9}$$

*and therefore any local maximum of $V_\pi$ is a stationary point of* (9).

(b) *If $\pi$ is a stationary point of* (9), *then $m_x + \frac{\tilde{\mu}}{2\pi_x}$ is constant for all $x \in E$, and $(\sum_{x \in J} m_x(\hat{\pi})\hat{\pi}_x)/(\sum_{x \in J} \hat{\pi}_x) + \tilde{\mu}/(2\sum_{x \in J} \hat{\pi}_x)$ is equal to the same constant for all $J \subset E$.*

(c) *Suppose $\hat{\pi}$ is a stationary point of* (9), *then $m_x(\hat{\pi}) \geq \overline{m}_{\hat{\pi}}$ (resp. $\leq$) if and only if $\hat{\pi}_x \geq 1/(2L+1)$ (resp. $\leq$). And $m_x(\hat{\pi}) > m_y(\hat{\pi})$ if and only if $\hat{\pi}_x > \hat{\pi}_y$.*

To say something specific about when speciation is likely to occur, we specialize to $m$, $K$ and $B$ of the following form in most of the results we establish (in fact, all results except Theorem 2.5):

ASSUMPTION 1. We have (1) $m_x$ is of the form defined in (4).
(2) $K : E \to (0, 1]$ symmetric and unimodal (i.e., increasing on $[-L, 0] \cap \mathbb{Z}$ and decreasing on $[0, L] \cap \mathbb{Z}$) with $K_0 = 1$.
(3) $B : \mathbb{Z} \to [0, 1]$ with $B_x = b + (1-b)1_{\{|x| \geq M\}}$ with $b \in [0, 1]$.

With the cooperation kernel as defined in 3 above, the individuals of phenotype $x$ are more efficient at using resource of types inside the interval $(x - M, x + M)$ than inside $[-L, x - M] \cup [x + M, L]$, and $b$ can be thought of as a measure of how efficiently individuals of phenotype $x$ use resources of types $(x - M, x + M)$, with $b = 0$ meaning maximally efficient.

If $b = 1$, then each individual can use resources of all types equally efficiently (or inefficiently), and every individual suffers exactly the same level of competition from the rest of the population. This actually means that competition plays no part in determining how fit site $x$ is and $m_x$ is proportional to $K_x$. Therefore, since $K_x$ is unimodal (hence $K_x$ is strictly increasing in $[-L, 0]$ and strictly decreasing in $[0, L]$), the fitness should be unimodal, too. Lemma 2.4(c) says that stationary distributions of (9) has the property of fitter sites having more mass, thus we expect the stationary distribution $\hat{\pi}$ to be unimodal as well. In particular, $\hat{\pi}$ should attain its maximum at $x = 0$. As $\tilde{\mu} \to 0$, Proposition 2.2 tells us that we can expect the peak of $\hat{\pi}$ concentrated around 0 to become sharper and sharper, approaching $\delta_0$, the $\delta$-measure concentrated at 0. In fact, Theorem 2.5 (where we do not assume Assumption 1) shows that if $B$ is of the form as in Assumption 1(3), then $b$ only needs to be close to 1 for this behavior to occur. In this case, speciation is not likely to occur.

THEOREM 2.5. *If $K : E \to (0, 1]$ is symmetric and unimodal with $K_0 = 1$, and $K_1 = K_{-1} < B_x \leq 1$ for all $x \in E$, then for any $\varepsilon > 0$, there exist $\mu$ and $1/N$ small enough such that $\nu^N(\{\pi \in \Delta : \pi_x > \varepsilon \text{ for some } x \in E \setminus \{0\}\}) < \varepsilon$, where $\nu^N$ is the stationary distribution of the Fleming–Viot model with generator $\mathcal{G}$.*

2.2. *Results for intense competition with relatively large $\mu$.* More interesting behavior arises when there is intense competition between pairs of sites that are close to each other, that is, $b$ is small and individuals of phenotype $x$ are far better at using resources of type $(x - M, x + M)$ than other types. In this case, speciation is likely to occur for certain $K$ and $B$. Whether or not speciation occurs depends on the exact shape of $K$ and the mutation parameter $\mu$, and is a difficult problem for general $K$ (even assuming symmetry and unimodality) and $\mu$.

We first present a result (Theorem 2.6 below) that roughly says that if any local maximum $\hat{\pi}$ of $V_\pi$ with $\hat{\pi}_0$ suitably small cannot have all the remaining mass on one side of 0, that is, there is significant amount of mass in both $[-L, -1]$ and $[1, L]$, which does not go to 0 as $\mu \to 0$.

THEOREM 2.6. *Suppose Assumption* 1 *holds. Let $M \geq 1$. If $\hat{\pi}$ is a stationary point of* (9) *where $\hat{\pi}_0 < K_1/(2M-1)(2L+1)$ and $\hat{\pi}_x > 1/(2L+1)$ for some site $x \in [-L, -1]$ (resp. $x \in [1, L]$), then*

$$\sum_{z=1}^{L} \hat{\pi}_z \geq \frac{K_1}{(2M-1)(2L+1)} \min\left(1, \frac{1}{2(1-b)}, \frac{1}{2(1-b)}\left(\frac{1}{K_1} - 1\right)\right)$$

*[resp. $\sum_{z=-L}^{-1} \hat{\pi}_z \geq K_1 \min(1/K_1 - 1, 1)/(2(1-b)(2M-1)(2L+1))$].*

REMARK 2.7. If $\hat{\pi}$ is a stationary point of (9) where $\hat{\pi}_0 < K_1/(2M-1) \times (2L+1) < 1/(2L+1)$, then some site other than 0 is bound to have more mass than the mean $1/(2L+1)$.

Let $\lfloor c \rfloor$ denote the largest integer less than or equal to $c$, and $\lceil c \rceil$ denote the smallest integer larger than or equal to $c$. Theorem 2.6 requires the mass at phenotype 0 to be rather small. Proposition 2.8 below guarantees this to be the case for any local maximum of $V_\pi$ if $\mu^{2/3}$ is relatively large compared to $b$ and the carrying capacity function $K_x$ decreases rapidly to very small levels (smaller than $\mu^{2/3}$) before reaching $M$ or $-M$.

PROPOSITION 2.8. *Define $l = -(n - M)$ and $p = \lceil M/2 \rceil$. Suppose Assumption* 1 *holds and $\hat{\pi}$ is a stationary point of* (9). *If $\tilde{\mu} \leq 4K_p^2/(4L+2)^3$ and there is an $n \leq M$ such that*

(10) $$b + K_n < (\tilde{\mu}K_p/4)^{2/3},$$

*then*

$$\hat{\pi}_x \leq \frac{\tilde{\mu}}{2((\tilde{\mu}K_p/4)^{2/3} - b - K_n)}$$

*for $x \in [-L, -n] \cup [-l, l] \cup [n, L]$.*

In the parameter regime of Proposition 2.8, the sites in $[-L,-n] \cup [n,L]$ have very small carrying capacity such that they cannot support significant population, and there is enough mutation effects to force most of the population into intervals $(-n,-l)$ and $(l,n)$ for any local maximum of $V_\pi$. This is a different effect from mutation simply spreading mass to all sites in $E$ evenly, since the mass in $(-n,-l) \cup (l,n)$ does not go to 0 as $\tilde{\mu} \to 0$, as long as $b+K_n < (\tilde{\mu}K_p/4)^{2/3}$ is satisfied. As shown by Proposition 2.10(b) below, $\delta_0$ is a local maximum of $\overline{m}_\pi$ if $K_M < b$, which may hold in the parameter regime of Proposition 2.8. But with the combined effect of mutation and selection, no configuration close to $\delta_0$ can be a local maximum of $V_\pi$. If for example $(b+K_n)/(\tilde{\mu}K_p/4)^{2/3} \to 0$, then Proposition 2.8 implies that for sufficiently small $\tilde{\mu}$, $\hat{\pi}_x \leq (16\tilde{\mu}/K_p^2)^{1/3}$ for $x \in [-L,-n] \cup [-l,l] \cup [n,L]$. Furthermore, the mass in $(-n,-l) \cup (l,n)$ cannot concentrate on one side of 0 by Theorem 2.6. Therefore in this case, speciation is likely to occur and there are at least 2 new species, with phenotypes in $(-n,-l)$ and $(l,n)$.

Proposition 2.8 holds even if $\delta_0$ is not a local maximum of $\overline{m}_\pi$, in which case Proposition 2.10(d) below provides possible existence (which is verified by our simulation) of a local maximum of $\overline{m}_\pi$ of form $\hat{\pi} = p\delta_{-M} + (1-2p)\delta_0 + p\delta_M$, where $p$ may be quite small. In this case, the relatively large mutation effects still prevent any configuration close to $\hat{\pi}$ from being a local maximum of $V_\pi$, and the population is driven toward a bimodal configuration.

2.3. *Results for intense competition with small $\mu$.* If $b$ is fixed but $\mu$ is sufficiently small, then it is possible for a configuration close to $\delta_0$ to be a local maximum of $V_\pi$. For such a parameter regime, we first present Theorem 2.9 below that says if there is little mass outside the intervals $(-\lfloor M/2 \rfloor, \lfloor M/2 \rfloor)$, then all mass must be concentrated at site 0 for sufficiently small $\tilde{\mu}$. This result precludes the existence of any configuration with mass spread amongst sites near 0 as a local maximum of $V_\pi$, such that if there is speciation, then there must be new phenotypes far away from 0 if $\tilde{\mu}$ is small.

THEOREM 2.9. *Let $q = \lfloor M/2 \rfloor$. If Assumption 1 holds and $\hat{\pi}$ is a stationary point of (9) that satisfies*

$$\sum_{x \in [-L,-q] \cup [q,L]} \hat{\pi}_x < \varepsilon < \frac{b(1-K_1)K_{q-1}}{b(1-K_1)K_{q-1}+K_1},$$

*then*

$$\sum_{x \in (-q,0) \cup (0,q)} \hat{\pi}_x \leq \frac{\tilde{\mu}}{2(b(1-K_1)K_{q-1}(1-\varepsilon)-K_1\varepsilon)}.$$

Simulations indicate that all local maxima of $\overline{m}_\pi$ have supports that consist of sites spread exactly $M$ apart, but a proof of this statement remains elusive. The landscape of $\overline{m}_\pi$ is rather complicated – explicit calculations for low dimensional systems ($L = 1$ and $L = 2$) indicate that there are interior stationary points of (6) but they are all saddle points of $\overline{m}_\pi$. Therefore if one tries to prove the above support property of local maxima, an approach that only looks at the stationary points of (6) would probably not work.

Proposition 2.10(a) below shows that the support of any local maximum $\hat{\pi}$ of $\overline{m}_\pi$ cannot have a support that consists of sites spread more than $M$ apart, and the rest of Proposition 2.10 gives formulas of local maxima with supports that consist of 1, 2 and 3 sites. Existence of local maxima with support that consists of more than 3 sites is also possible, but then explicit calculation becomes prohibitive. If $\tilde{\mu}$ is sufficiently small, Proposition 2.3 provides local maxima of $V_\pi$ that are close to the local maxima of $\overline{m}_\pi$.

PROPOSITION 2.10. *Suppose Assumption 1 holds and $\hat{\pi}$ is a stationary point of* (9).

(a) *If $\hat{\pi}$ has mass only on $x_1 < \cdots < x_k$ and $x_{j+1} - x_j \geq M + 1$ for some $j$, then $\hat{\pi}$ cannot be a local maximum of $\overline{m}_\pi$.*

(b) *If $K_M = K_{-M} < b$, then $\delta_0$ is the unique local maximum of $\overline{m}_\pi$ in a small neighborhood of $\delta_0$. If $K_M = K_{-M} > b$, then $\delta_0$ is not a local maximum of $\overline{m}_\pi$.*

(c) *Let $x \in [-M+1, -1] \cap \mathbb{Z}$ and $p = K_{-x+M}(K_{-x} - bK_{-x+M})/(2K_{-x} \times K_{-x+M} - bK_{-x}^2 - bK_{-x+M}^2)$, then $\hat{\pi} = p\delta_{-x} + (1-p)\delta_{-x+M}$ is the unique local maximum of $\overline{m}_\pi$ in a small neighborhood of $\hat{\pi}$ with*

$$\overline{m}_{\hat{\pi}} = \frac{K_{-x}^2 K_{-x+M}^2 (1-b^2)}{2K_{-x}K_{-x+M} - bK_{-x}^2 - bK_{-x+M}^2},$$

*if $b$, $K_{-x-M}$, and $K_{-x+2M}$ are all $< K_{-x}K_{-x+M}(1+b)/(K_{-x} + K_{-x+M})$, and $bK_{-x}^2 + bK_{-x+M}^2 < 2K_{-x}K_{-x+M}$.*

(d) *Let $x \in [-M+1, 0] \cap \mathbb{Z}$ (resp., $x \in (0, M-1] \cap \mathbb{Z}$),*

$$a = K_{x-M}^2 K_x^2 + K_{x-M}^2 K_{x+M}^3 + K_x^2 K_{x+M}^3,$$

$$c = 2K_{x-M}K_x K_{x+M}(K_{x-M} + K_x + K_{x+M}) - (1+b)a,$$

$$p = \frac{K_x K_{x+M}}{c}(K_{x-M}K_{x+M} + K_{x-M}K_x - (1+b)K_x K_{x+M}),$$

$$q = \frac{K_{x-M}K_{x+M}}{c}(K_{x-M}K_x + K_x K_{x+M} - (1+b)K_{x-M}K_{x+M}),$$

*then $\hat{\pi} = p\delta_{x-M} + q\delta_x + (1-p-q)\delta_{x+M}$ is the unique local maximum of $\overline{m}_\pi$ in a small neighborhood of $\hat{\pi}$ with*

$$\overline{m}_{\hat{\pi}} = \frac{(2-b-b^2)K_{x-M}^2 K_x^2 K_{x+M}^2}{2K_{x-M}K_x K_{x+M}(K_{x-M} + K_x + K_{x+M}) - (1+b)a},$$

*if*

$$2 - b > (1 - b)\left(\frac{1}{K_x} + \frac{1}{K_{x-M}}\right) - \frac{1}{K_{x+M}}$$

$$\left[resp., \ 2 - b > (1 - b)\left(\frac{1}{K_x} + \frac{1}{K_{x+M}}\right) - \frac{1}{K_{x-M}}\right],$$

$$K_{x-2M} < \frac{K_{x-M} K_x K_{x+M}(2 - b - b^2)}{(1-b)(K_{x-M}K_x + K_{x-M}K_{x+M} + K_x K_{x+M})},$$

$$K_{x+2M} < \frac{K_{x-M} K_x K_{x+M}(2 - b - b^2)}{(1-b)(K_{x-M}K_x + K_{x-M}K_{x+M} + K_x K_{x+M})}$$

*and*

$$1 + b < \frac{2}{a} K_{x-M} K_x K_{x+M}(K_{x-M} + K_x + K_{x+M}).$$

2.4. *Simulation and discussion.* Even though we have no general recipe for identifying the global maximum of $V_\pi$, but only a class of local maxima, it is nevertheless possible to check whether speciation is likely to occur. For that, we can first check whether a configuration close to $\delta_0$ is a local maximum of $V_\pi$; if not, then speciation is likely to occur.

As $\tilde{\mu} \to 0$, local maxima of $V_\pi$ converge to those of $m_\pi$. If we are interested in the behavior of the process when $\tilde{\mu}$ is in this regime, it suffices to check whether $\delta_0$ is a local maximum of $m_\pi$, and Proposition 2.10(b) provides the answer.

If we are interested in the behavior of the process when $\tilde{\mu}$ is small but fixed, then Theorem 2.6 and Theorem 2.9 provide partial answers. If there is no configuration close to $\delta_0$ which is a local maximum of $V_\pi$, then by Theorem 2.9 (its contrapositive), there must be significant mass in $[-L, -\lfloor M/2 \rfloor] \cup [\lfloor M/2 \rfloor, L]$ for sufficiently small $\tilde{\mu}$, so that there are new phenotypes far away from 0. Moreover, if $\hat{\pi}_0$ is quite small ($< K_1/(2M - 1)(2L + 1)$), then Theorem 2.6 says that there must be significant mass on both sides of site 0. These two results taken together mean that if configurations close to $\delta_0$ is not a local maximum of $V_\pi$, then for any local maximum, there is significant (i.e., does not get small when $\tilde{\mu}$ gets very small) mass at more than one site, hence speciation is likely to occur.

Proposition 2.8 provides a condition under which all local maxima of $V_\pi$ are bimodal. In this case, $\mu^{2/3}$ is relatively large compared to $b$ and the drift component (terms involving first derivatives) of $\mathcal{G}$ drives the population toward a bimodal configuration. An example of this case is shown in Figure 4.

On the other hand, if $\mu$ is small and $\delta_0$ is a local maximum of $\overline{m}_\pi$, then a configuration close to $\delta_0$ is a local maximum of $V_\pi$. But it does not mean that speciation is unlikely to occur since $\delta_0$ may not be the global maximum of

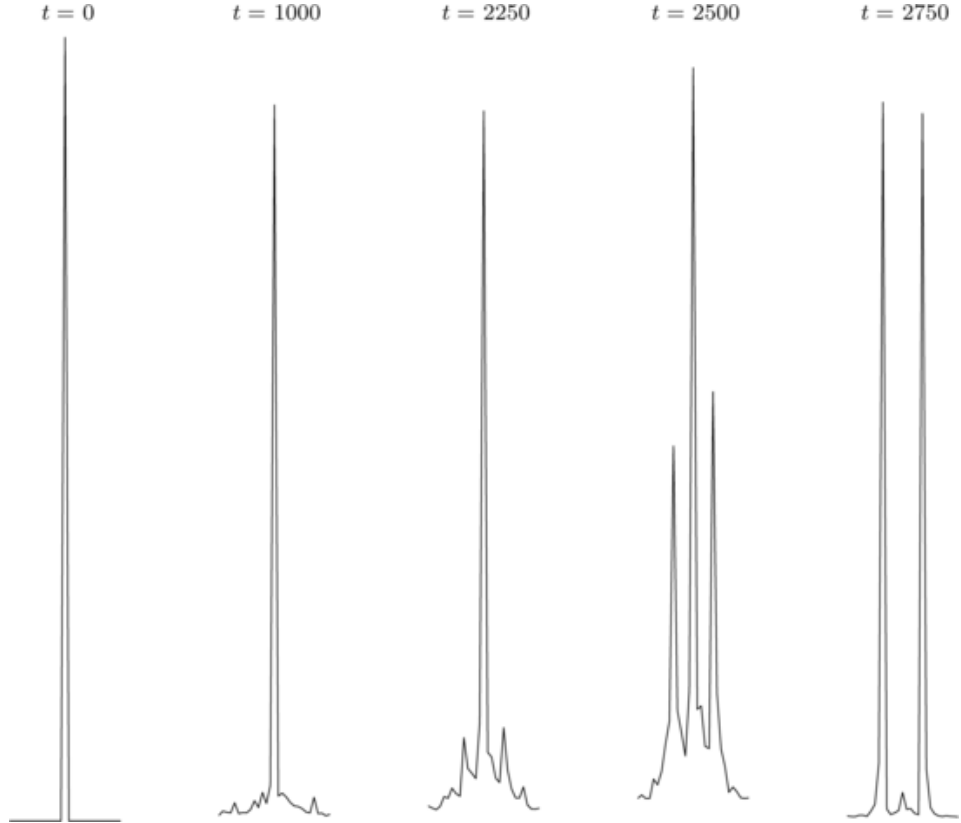

FIG. 4. *When $\mu$ is relatively large, $\delta_0$ is not a local maximum and speciation occurs at time $t = 2750$. Here $L = 14$, $N = 225^2$, $B_x = 0.01 + 0.99 \cdot 1_{\{|x| \geq 10\}}$, $K_x = \exp(-x^2/20)$ and $\mu = 6 \times 10^{-5}$.*

$\overline{m}_\pi$, Proposition 2.10(c) and (d) provides other local maxima of $m_\pi$. If one of these local maxima, not $\delta_0$, is the global maximum of $\overline{m}_\pi$, then the stationary distribution $\nu^N$ concentrates mostly on a configuration that has at least 2 modes, and speciation is likely to occur. But if the initial configuration is very close to $\delta_0$, then the drift component of $\mathcal{G}$ moves $\pi$ toward the local maximum of $V_\pi$ that resembles $\delta_0$. One must wait long enough for the noise component (the term involving second derivatives) of $\mathcal{G}$ to drive $\pi$ away from $\delta_0$ and to a configuration in the basin of attraction of the more fit bimodal configuration. Since the noise component of $\mathcal{G}$ is $O(1/N)$, it may take a long time to get away from the less fit local maximum that resembles $\delta_0$ when the population size is large, that is, speciation may take much longer to occur than in the case of $\delta_0$ not being a local maximum of $\overline{m}_\pi$. An example of this case is shown in Figure 5.

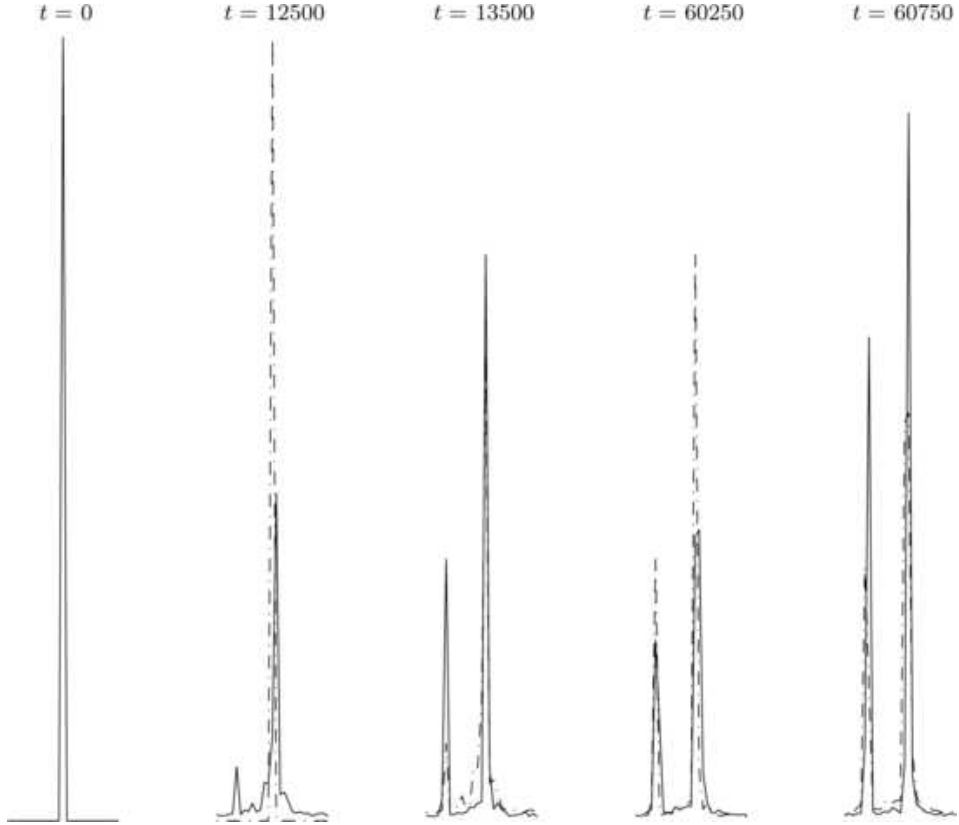

FIG. 5. *When $\mu$ is small, $\delta_0$ is a local maximum and speciation occurs at $t = 13500$, much later than in Figure 4. Even after speciation occurs, the population continues to move to fitter local maxima, but this takes even longer. Here $L = 14$, $N = 225^2$, $B_x = 0.01 + 0.99 \cdot 1_{\{|x| \geq 10\}}$, $K_x = \exp(-x^2/20)$ and $\mu = 5 \times 10^{-5}$. For comparison purposes, the dotted lines at $t = 12{,}500$ to $t = 60750$ are duplicates of the figures on their left.*

As can be seen from the figure above, speciation occurs around time $t = 13500$, but this is a different configuration from the one reached in Figure 4 ($\mu = 6 \times 10^{-5}$) at time $t = 2750$. The configuration in Figure 4 at time $t = 2750$ has the two peaks symmetrically placed on both sides of $x = 0$, at $x = 5$ and $x = -5$, but the configuration in Figure 5 ($\mu = 5 \times 10^{-5}$) at time $t = 13{,}500$ has the one of the peaks at $x = 1$ and the other at $x = -9$, a much less fit configuration than the one with symmetrically placed peaks. At time $t = 60{,}750$, the selection component of $\mathcal{G}$ succeeds in driving $\pi$ to a fitter configuration, with two peaks placed at $x = 2$ and $x = -8$, but with $\mu = 5 \times 10^{-5}$ instead of $\mu = 6 \times 10^{-5}$, it will take much longer to reach a configuration with peaks at $x = -5$ and $x = 5$, supposedly the global maximum.

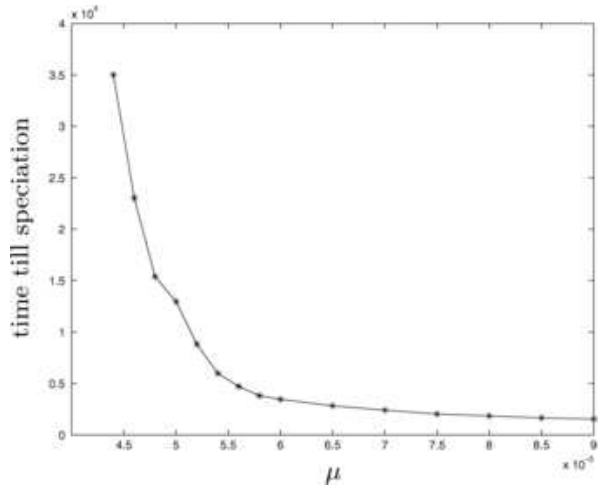

Fig. 6. *A plot of speciation time (mean of* 20–50 *realizations) against the mutation parameter* $\mu$. *Here* $L = 14$, $N = 225^2$, $B_x = 0.01 + 0.99 \cdot 1_{\{|x| \geq 10\}}$ *and* $K_x = \exp(-x^2/20)$. *We did not perform simulations for* $\mu < 4.4 \times 10^{-5}$ *because it takes too long, but we expect speciation time to grow at a fast rate as* $\mu$ *decreases beyond* $4.4 \times 10^{-5}$.

Simulations indicate that for $\mu > 5.028$ and $K$, $B$ and $N$ as given in Figures 4 and 5, $\delta_0$ is a local maximum of $V_\pi$, but $\delta_0$ is not a local maximum of $V_\pi$ if $\mu < 5.027$. We surmise that the deterministic dynamical system (9) has a bifurcation near $\mu = \check{\mu}$ ($\tilde{\mu} \approx 5.027$ for the simulations shown in Figures 4 and 5), which causes the drastically different speciation time of the Fleming–Viot process (5) when $\mu$ decreases from $5.5 \times 10^{-5}$ to $4.5 \times 10^{-5}$ (see Figure 6). Notice that our simulation has $N = 225^2$, so that the noise component is large enough for the simulation to achieve speciation in reasonable amount of time when $\mu < \check{\mu}$. When $N$ is much larger than $225^2$, we expect the increase in the time until speciation to be much more drastic. And if $N = \infty$, then speciation will never occur for $\mu < \check{\mu}$ if one starts with initial condition $\delta_0$, since in this case, the dynamical system (9) coincides with (7), the infinite population limit of the Fleming–Viot process.

Thus whether or not a configuration close to $\delta_0$ is local maximum of $V_\pi$ affects the length of time it takes for speciation to occur, assuming that a configuration with significant mass at sites other than 0 is the global maximum of $V_\pi$ such that speciation is likely to occur. In addition to the unexpected dependence on the mutation parameter $\mu$, the carrying capacity $K$ and cooperation kernel $B$ also affects greatly whether a certain configuration is a local maximum of $\overline{m}_\pi$ or $V_\pi$. In Proposition 2.10(b), for example, a small change in $K$ or $b$ can change whether $\delta_0$ is a local maximum of

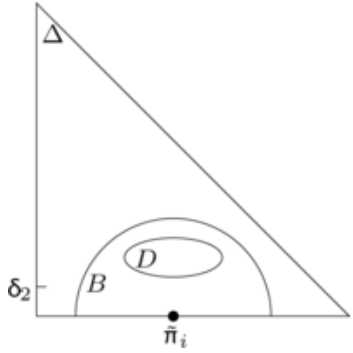

Fig. 7. *How to pick D.*

$\overline{m}_\pi$. This is somewhat similar to the behavior exhibited by the conditioned Dieckmann–Doebeli model, where we remarked near Figures 2 and 3 that a slight change in $K$ or $C$ may have a large effect on the shape of stationary distributions of (1).

**3. Proofs.** We start by proving Proposition 2.2 and Proposition 2.3, the two statements that relate global/local maximum of $V_\pi$ to that of $m_\pi$, when $\mu$ is small.

PROOF OF PROPOSITION 2.2. We observe that $\overline{m} : \pi \mapsto \overline{m}_\pi$ is a *quadratic* function, hence a continuous and open mapping. Near each global maximum $\tilde{\pi}_i$, the Hessian matrix must be positive definite; otherwise, there would exist an entire subspace of global maxima. Therefore for sufficiently small $\varepsilon$, we can pick open neighborhoods $A$ and $B$ such that $B \subset A \subset \bigcup_{i=1}^{k} \text{Ball}(\tilde{\pi}_i, \varepsilon)$ and

$$\inf_{\pi \in B} \overline{m}_\pi - \sup_{\pi \in \Delta \setminus A} \overline{m}_\pi > \delta$$

for some positive $\delta$. Since the $\pi \in \Delta$ that maximizes $\sum_{x=-L}^{L} \log \pi_x$ places equal weights on $x \in E$,

(11) $$\sup_{x \in \Delta} \sum_{x=-L}^{L} \log \pi_x \leq -(2L+1)\log(2L+1).$$

Furthermore, there exist positive $\delta_1$ and $\delta_2$ (independent of $\mu$ and $N$) and an open set $D \subset B$, even if some $\tilde{\pi}_i$'s are on the boundary of $\Delta$, such that $|D| > \delta_1$ (where $|D|$ is the volume of $D$) and $\pi_x > \delta_2$ for all $\pi \in D$ and $x \in E$ (see Figure 7), thus

$$\inf_{x \in D} \sum_{x=-L}^{L} \log \pi_x \geq -(2L+1)\log \frac{1}{\delta_2}.$$

Combining the three inequalities above, we obtain

$$\inf_{\pi \in D} V_\pi - \sup_{\pi \in \Delta \setminus A} V_\pi > \delta + \tilde{\mu}(2L+1)\left(-\log \frac{1}{\delta_2} + \log(2L+1)\right).$$

Consequently, if $\mu$ and $1/N$ are sufficiently small, then

$$\inf_{\pi \in D} V_\pi - \sup_{\pi \in \Delta \setminus A} V_\pi > \frac{\delta}{2},$$

and thus

$$\frac{\inf_{\pi \in D} \exp((N/2)V_\pi)}{\sup_{\pi \in \Delta \setminus A} \exp((N/2)V_\pi)} > e^{N\delta/4}.$$

Define $C_1 = \sup_{\pi \in \Delta \setminus A} \exp(NV_\pi/2)$, then we obtain the following bounds for $\int \exp(NV_\pi/2) \, d\pi$:

$$\int_A \exp\left(\frac{N}{2}V_\pi\right) d\pi > \int_D \exp\left(\frac{N}{2}V_\pi\right) d\pi > \delta_1 C_1 e^{N\delta/4},$$

$$\int_{\Delta \setminus A} \exp\left(\frac{N}{2}V_\pi\right) d\pi < C_1 |\Delta|,$$

where we use $|D| > \delta_1$ to obtain the first inequality above. The two inequalities above in turn imply that

$$\nu^N(A) = \frac{\int_A \exp((N/2)V_\pi) \, d\pi}{\int_\Delta \exp((N/2)V_\pi) \, d\pi} \geq \frac{\delta_1 e^{N\delta/4}}{|\Delta| + \delta_1 e^{N\delta/4}}$$

$$= 1 - \frac{|\Delta|}{|\Delta| + \delta_1 e^{N\delta/4}} > 1 - \varepsilon,$$

if $N$ is sufficiently large. Since $A \subset \bigcup_{i=1}^k \text{Ball}(\tilde{\pi}_i, \varepsilon)$, we are done. □

Part of the proof above can be easily adapted to prove a related statement on local maxima.

PROOF OF PROPOSITION 2.3. We use similar ideas as in the proof of Proposition 2.2. Near $\tilde{\pi}$, the Hessian matrix of the quadratic function $\overline{m}_\pi$ must be positive definite. Therefore we can pick open neighborhoods $A$, $B$ and $Z$ such that $B \subset A \subset Z$ of $\tilde{\pi}$ and

$$\inf_{\pi \in B} \overline{m}_\pi - \sup_{\pi \in Z \setminus A} \overline{m}_\pi > \delta$$

for some positive $\delta$. Furthermore, there exists positive $\delta_2$ (independent of $\mu$), and an open set $D \subset B$ such that $\pi_x > \delta_2$ for all $\pi \in D$ and $x \in E$, and thus

$$\inf_{x \in D} \sum_{x=-L}^{L} \log \pi_x \geq -(2L+1) \log \frac{1}{\delta_2}.$$

The two inequalities above and (11) imply

$$\inf_{\pi\in D} V_\pi - \sup_{\pi\in Z\setminus A} V_\pi > \delta + \tilde{\mu}(2L+1)\left(-\log\frac{1}{\delta_2} + \log(2L+1)\right).$$

Consequently, if $\tilde{\mu}$ is sufficiently small, then

$$\inf_{\pi\in D} V_\pi - \sup_{\pi\in Z\setminus A} V_\pi > \frac{\delta}{2},$$

and there exists a local maximum of $V_\pi$ in $A$. $\square$

We study local maxima of $V_\pi$ by checking various points in $\Delta$ for stationarity when evolved according to (9). It also happens that the dynamical system (9) is almost the same as the limiting dynamical system of our Fleming–Viot process (7).

PROOF OF LEMMA 2.4. (a) Notice that $\partial_{\pi_x}\overline{m}_\pi = 2\sum_z K_x B_{x-z} K_z \pi_z = 2m_x$. Therefore if $\pi$ evolves according to (9), then

$$\partial_t V_\pi = \sum_x (\partial_{\pi_x} V_\pi)(\partial_t \pi_x)$$

$$= \sum_x \left(2m_x + \frac{\tilde{\mu}}{\pi_x}\right)\pi_x\left(m_x + \frac{\tilde{\mu}}{2\pi_x} - \overline{m}_\pi - \frac{\tilde{\mu}}{2}(2L+1)\right)$$

$$= 2\sum_x \left(m_x + \frac{\tilde{\mu}}{2\pi_x}\right)\pi_x\left(m_x + \frac{\tilde{\mu}}{2\pi_x} - \overline{m}_\pi - \frac{\tilde{\mu}}{2}(2L+1)\right)$$

$$- \left(\overline{m}_\pi + \frac{\tilde{\mu}}{2}(2L+1)\right)\sum_x \pi_x\left(m_x + \frac{\tilde{\mu}}{2\pi_x} - \overline{m}_\pi - \frac{\tilde{\mu}}{2}(2L+1)\right),$$

where the term in the last line $\sum_x \pi_x(m_x + \frac{\tilde{\mu}}{2\pi_x} - \overline{m}_\pi - \frac{\tilde{\mu}}{2}(2L+1)) = \sum_x \pi_x \times m_x - \overline{m}_\pi$, which is zero. This implies

$$\partial_t V_\pi = 2\sum_x \pi_x\left(m_x + \frac{\tilde{\mu}}{2\pi_x} - \overline{m}_\pi - \frac{\tilde{\mu}}{2}(2L+1)\right)^2 \geq 0,$$

which establishes (a).

(b) From (a), at any stationary point $\hat{\pi}$ of (9), we have $m_x(\hat{\pi}) - \overline{m}_{\hat{\pi}} + \frac{\tilde{\mu}}{2}(1/\hat{\pi}_x - (2L+1)) = 0$ for all $x \in E$, therefore $m_x(\hat{\pi}) + \tilde{\mu}/2\hat{\pi}_x$ is constant for all $x \in E$. We define this constant to be $c = m_x(\hat{\pi}) + \tilde{\mu}/2\hat{\pi}_x$, then $m_x(\hat{\pi})\hat{\pi}_x + \tilde{\mu}/2 = c\hat{\pi}_x$ for all $x \in E$. Thus for all $J \subset E$,

$$\frac{\sum_{x\in J} m_x(\hat{\pi})\hat{\pi}_x}{\sum_{x\in J} \hat{\pi}_x} + \frac{\tilde{\mu}}{2\sum_{x\in J} \hat{\pi}_x} = c.$$

(c) Since

$$m_x(\hat{\pi}) - \overline{m}_{\hat{\pi}} + \frac{\tilde{\mu}}{2}\left(\frac{1}{\hat{\pi}_x} - (2L+1)\right) = 0, \tag{12}$$

$m_x(\hat{\pi}) \geq \overline{m}_{\hat{\pi}}$ implies $1/\hat{\pi}_x - (2L+1) \leq 0$ and vice versa. The cases of $m_x(\hat{\pi}) \leq \overline{m}_{\hat{\pi}}$ and $m_x(\hat{\pi}) > m_y(\hat{\pi})$ are similar. $\square$

In fact, according to Theorem A.9 of [1], (9) is a so-called Svirezhev–Shahshahani gradient system with potential $V_\pi$, that is, $\partial_t \pi = \tilde{\nabla} V(\pi)$, where $\tilde{\nabla} V(\pi) = G_\pi \nabla V(\pi)$ and $G_\pi$ is the matrix formed by entries $g^{xy} = \pi_x(\delta_{xy} - \pi_y)$. Any gradient system, such as (9), has the property that all orbits, regardless of initial condition, converge to some point in the $\omega$-limit set

$$D_\omega = \{p : p \text{ is an accumulation point of } \pi(t) \text{ as } t \to \infty\}.$$

All points in $D_\omega$ are stationary points of (9).

3.1. *Mild competition.* In this section, we establish Theorem 2.5, which says that speciation is impossible if competition between phenotypes close to each other is mild.

LEMMA 3.1. *If $\hat{\pi}$ is a local maximum of $\overline{m}_\pi$ with support $S$, that is, $\hat{\pi}_x > 0$ for $x \in S$ only, then $m_x(\hat{\pi})$ are all equal for $x \in S$. In other words, if $m_x(\hat{\pi}) \neq m_y(\hat{\pi})$, then either $\hat{\pi}_x = 0$ or $\hat{\pi}_y = 0$.*

PROOF. A simple calculation involving Lagrange multipliers:

$$\frac{\partial}{\partial \pi_x}\left(\overline{m}_\pi + \lambda \sum_{x \in S} \pi_x\right) = 0 \implies 2m_x + \lambda = 0$$

establishes the desired result. $\square$

This observation enables us to establish the following:

LEMMA 3.2. *If $K : E \to (0, 1]$ is symmetric and unimodal with $K_0 = 1$, and $K_1 = K_{-1} < B_x \leq 1$ for all $x \in E$, then $\delta_0$ is the unique global maximum of the mean fitness function $\overline{m}_\pi$.*

PROOF. Assume $K_1 = K_{-1} < B_x \leq 1$ for all $x \in E$. The following two estimates

$$m_0 = K_0 \sum_z B_{-z} K_z \pi_z \geq K_0 \left(\inf_z B_z\right) \sum_z K_z \pi_z = \left(\inf_z B_z\right) \sum_z K_z \pi_z,$$

$$m_x = K_x \sum_z B_{x-z} K_z \pi_z \leq K_x \sum_z K_z \pi_z \quad \text{for } x \neq 0$$

imply that $m_0 > m_x$ for all $x \neq 0$ if $\inf_z B_z > \sup_{x \neq 0} K_x = K_1$. Lemma 3.1 then implies that any local maximum $\hat{\pi}$ of $\overline{m}_\pi$ whose support includes 0 must have either $\hat{\pi}_0 = 1$ or $\hat{\pi}_0 = 0$. If $\hat{\pi}_0 = 1$, then $\hat{\pi} = \delta_0$, and $\overline{m}_{\hat{\pi}} = m_0(\hat{\pi}) = B_0 > K_1$. But if $\hat{\pi}_0 = 0$, then for all $x$ in the support of $\hat{\pi}$,

$$\overline{m}_{\hat{\pi}} = m_x \leq K_1 \sum_z K_z \pi_z \leq K_1 < B_0$$

by assumption. Therefore a local maximum $\hat{\pi}$ whose support does not include 0 must have smaller mean fitness than $\delta_0$, that is, $\delta_0$ is the unique global maximum of $\overline{m}_\pi$. □

PROOF OF THEOREM 2.5. A direct application of Proposition 2.2 leads to the desired result, which actually applies to very general $B$ (but still requires $K$ to be symmetric and unimodal). □

3.2. *Intense competition.* Now we focus on the case of intense competition, that is, $b$ is small. The goal is to establish some conditions under which speciation is likely to occur. The first result in this direction is Theorem 2.6, which roughly says that any local maximum $\hat{\pi}$ of $V_\pi$ with $\hat{\pi}_0$ suitably small cannot have all remaining mass on one side of 0. For this, we first establish the following lemma, which assumes that there is significant mass at site $-y$ to the left of site 0. If $-y > -M$, then the lemma shows there is mass in the interval $[-y+M, M]$, which implies Theorem 2.6. But if $-y \leq -M$, then the lemma shows there is mass in the interval $[-y+M, 2-M]$, which implies that the previous case holds, and in turn Theorem 2.6 holds as well.

LEMMA 3.3. *Let $M \geq 1$. Suppose Assumption 1 holds and $\hat{\pi}$ is a stationary point of (9).*

(a) *If $\hat{\pi}_{-y} \geq 1/(2L+1)$ for some $-y \in [-L, -M]$ and $\hat{\pi}_z < 1/(2L+1)$ for all $z \in [-y+1, -1]$, then $\sum_{z=2-M}^{M} \hat{\pi}_z \geq \sum_{z=\max(-y+M,2-M)}^{M} \hat{\pi}_z \geq K_1/(2L+1)$.*

(b) *If $\hat{\pi}_{-y} \geq D \geq \hat{\pi}_0$ for some $-y \in [1-M, -1]$, then*

$$\sum_{z=-y+M}^{M-1} \hat{\pi}_z \geq \min\left(D, \frac{D}{2(1-b)}, \frac{D}{2(1-b)}\left(\frac{1}{K_1} - 1\right)\right).$$

PROOF. (a) We observe that $-y + M < 1$, therefore the cooperation intensity between sites $-y$ and 1 is 1. We define

$$A = \sum_{z=-L}^{-y-M} K_z \hat{\pi}_z + b \sum_{z=-y-M+1}^{M} K_z \hat{\pi}_z + \sum_{z=M+1}^{L} K_z \hat{\pi}_z + (1-b) \sum_{z=-y+M}^{1-M} K_z \hat{\pi}_z,$$

where the last sum may be over an empty set, in which case that sum is defined to be 0. The fitness of sites $-y$ and $1$ are

$$m_{-y}(\hat{\pi}) = K_{-y}\left(A + (1-b)\sum_{z=\max(-y+M,2-M)}^{M} K_z\hat{\pi}_z\right),$$

$$m_1(\hat{\pi}) = K_1\left(A + (1-b)\sum_{z=-y-M+1}^{\min(1-M,-y+M-1)} K_z\hat{\pi}_z\right).$$

If $m_{-y}(\hat{\pi}) < m_1(\hat{\pi})$, then $\hat{\pi}_1 > \hat{\pi}_{-y} \geq 1/(2L+1) > K_1/(2L+1)$ by Lemma 2.4(c) and we are done. Otherwise, $m_{-y}(\hat{\pi}) \geq m_1(\hat{\pi})$ and since $K_{-y} \leq K_1$, the two equations above imply

$$K_{-y}\sum_{z=\max(-y+M,2-M)}^{M} K_z\hat{\pi}_z \geq K_1 \sum_{z=-y-M+1}^{\min(1-M,-y+M-1)} K_z\hat{\pi}_z \geq \frac{K_1 K_{-y}}{2L+1},$$

therefore

$$\sum_{z=\max(-y+M,2-M)}^{M} \hat{\pi}_z \geq \sum_{z=\max(-y+M,2-M)}^{M} K_z\hat{\pi}_z \geq \frac{K_1}{2L+1},$$

as required.

(b) We define

$$B = \sum_{z=-L}^{-y-M} K_z\hat{\pi}_z + b\sum_{z=-y-M+1}^{M-1} K_z\hat{\pi}_z + \sum_{z=M}^{L} K_z\hat{\pi}_z,$$

then the fitness of sites $-y$ and $0$ (the cooperation intensity between these sites is $b$) are:

(13) $$m_{-y}(\hat{\pi}) = K_{-y}\left(B + (1-b)\sum_{z=-y+M}^{M-1} K_z\hat{\pi}_z\right),$$

(14) $$m_0(\hat{\pi}) = K_0\left(B + (1-b)\sum_{z=-y-M+1}^{-M} K_z\hat{\pi}_z\right).$$

If $m_{-y}(\hat{\pi}) < m_{-y+M}(\hat{\pi})$, then $\hat{\pi}_{-y+M} > \hat{\pi}_{-y}$ by Lemma 2.4(c), and we are done. Otherwise,

$$m_{-y}(\hat{\pi}) \geq m_{-y+M}(\hat{\pi}) \geq K_{-y+M}K_{-y}\hat{\pi}_{-y} \geq K_{-y+M}K_{-y}D.$$

The above inequality and (13) imply that either

$$K_{-y}B \geq \frac{K_{-y+M}K_{-y}D}{2}$$

or
$$K_{-y}(1-b)\sum_{z=-y+M}^{M-1} K_z \hat{\pi}_z \geq \frac{K_{-y+M}K_{-y}D}{2},$$

or both. If $K_{-y}(1-b)\sum_{z=-y+M}^{M-1} K_z\hat{\pi}_z \geq K_{-y+M}K_{-y}D/2$, then

$$\frac{K_{-y+M}D}{2(1-b)} \leq \sum_{z=-y+M}^{M-1} K_z \hat{\pi}_z \leq K_{-y+M}\sum_{z=-y+M}^{M-1} \hat{\pi}_z$$

since $M-1 \geq -y+M > 0$, hence

(15) $$\sum_{z=-y+M}^{M-1} \hat{\pi}_z \geq \frac{D}{2(1-b)}.$$

And if $K_{-y}B \geq K_{-y+M}K_{-y}D/2$, then

(16) $$B \geq \frac{K_{-y+M}D}{2}.$$

Since $\hat{\pi}_0 \leq D \leq \hat{\pi}_{-y}$, we have $m_{-y}(\hat{\pi}) \geq m_0(\hat{\pi})$, then (13) and (14) imply the following:

$$1 \leq \frac{m_{-y}(\hat{\pi})}{m_0(\hat{\pi})} = \frac{K_{-y}(B+(1-b)\sum_{z=-y+M}^{M-1} K_z\hat{\pi}_z)}{K_0(B+(1-b)\sum_{z=-y-M+1}^{-M} K_z\hat{\pi}_z)},$$

$$\frac{1}{K_{-y}} \leq \frac{B+(1-b)\sum_{z=-y+M}^{M-1} K_z\hat{\pi}_z}{B} = 1 + \frac{1-b}{B}\sum_{z=-y+M}^{M-1} K_z\hat{\pi}_z,$$

$$\sum_{z=-y+M}^{M-1} K_z\hat{\pi}_z \geq \frac{B(1-K_{-y})}{(1-b)K_{-y}} \geq \frac{(1-K_{-y})K_{-y+M}D}{2(1-b)K_{-y}}$$

by (16). Since $\max_{z\in[-y+M,M-1]} K_z = K_{-y+M}$, we obtain

(17) $$\sum_{z=-y+M}^{M-1} \hat{\pi}_z \geq \frac{(1-K_{-y})D}{2(1-b)K_{-y}}.$$

Inequalities (15) and (17) imply the desired result. $\square$

PROOF OF THEOREM 2.6. We prove this for the case of $x \in [-L,-1]$. The other case of $x \in [1,L]$ is similar. Let $-y$ be the right most (i.e., the largest) site in $[-L,-1]$ where $\pi_{-y} > 1/(2L+1)$. We distinguish two cases, $-y+M < 1$ and $-y+M \geq 1$.

If $-y+M < 1$, then Lemma 3.3(a) implies that at least one site $y' \in [2-M,M]$ has more mass than $K_1/(2M-1)(2L+1)$. If $y' \in [1,M]$ then we are done; $y'$ cannot be 0 by assumption; and if $y' \in [2-M, -1]$, then we use the next paragraph to establish the result.

Now we define $D = K_1/(2M-1)(2L+1)$ and deal with the case of $\hat{\pi}_{-y} \geq D$ for some $-y + M \geq 1$. This includes the case $\pi_{-y} > 1/(2L+1)$ where $-y + M \geq 1$. We only need to apply Lemma 3.3(b) to reach the desired conclusion. $\square$

3.2.1. *Relatively large $\mu$.* For $\mu$ suitably large compared to $b$, we establish Proposition 2.8, which says that most of the mass is forced into 2 intervals on both sides of site 0, with little mass everywhere else. Theorem 2.6 implies that there must be significant mass in both these intervals, therefore speciation is likely to occur. For this, we first establish a lemma that gives a lower bound on the mean fitness $\overline{m}_\pi$ in terms of the mutation parameter $\tilde{\mu}$.

LEMMA 3.4. *Let $p = \lceil M/2 \rceil$. Suppose Assumption 1 holds and $\hat{\pi}$ is a stationary point of (9). If $\tilde{\mu} \leq 4K_p^2/(4L+2)^3$, then $\overline{m}_{\hat{\pi}} \geq (\tilde{\mu} K_p/4)^{2/3}$.*

PROOF. If $\hat{\pi}$ is a stationary point of (9), then since $m_x \geq 0$ for any $x$, (12) implies
$$\frac{1}{\hat{\pi}_x} - (2L+1) \leq \frac{2\overline{m}_{\hat{\pi}}}{\tilde{\mu}},$$
hence
$$\hat{\pi}_x \geq \frac{1}{(2\overline{m}_{\hat{\pi}}/\tilde{\mu}) + 2L + 1}.$$
Since $2p \geq M$, we have $B_{p-(-p)} = 1$, therefore
$$\overline{m}_{\hat{\pi}} \geq \hat{\pi}_{-p} K_p^2 \hat{\pi}_p \geq \left(\frac{K_p}{(2\overline{m}_{\hat{\pi}}/\tilde{\mu}) + 2L + 1}\right)^2 \geq \min\left(\frac{\tilde{\mu} K_p}{4\overline{m}_{\hat{\pi}}}, \frac{K_p}{4L+2}\right)^2$$
and
$$\overline{m}_{\hat{\pi}} \geq \min\left(\left(\frac{\tilde{\mu} K_p}{4}\right)^{2/3}, \left(\frac{K_p}{4L+2}\right)^2\right).$$
If $\tilde{\mu} \leq 4K_p^2/(4L+2)^3$, then the above estimate reduces to the desired conclusion. $\square$

PROOF OF PROPOSITION 2.8. We estimate the fitness of sites near 0. For $x \in [-l, l]$,
$$m_x(\hat{\pi}) = bK_x \sum_z K_z \hat{\pi}_z + (1-b)K_x \sum_{z \notin [x-M+1, x+M-1]} K_z \hat{\pi}_z$$

$$
\begin{align}
&\leq bK_0 \sum_z K_0 \hat{\pi}_z + K_0 K_n \sum_{z \notin [x-M+1, x+M-1]} \hat{\pi}_z \tag{18} \\
&\leq b + K_n, \tag{19}
\end{align}
$$

where we use the fact that $K_x$ is decreasing in $[0, L]$ in the second line. And for $x \notin [-n+1, n-1]$,

$$m_x(\hat{\pi}) \leq K_x \leq K_n. \tag{20}$$

Define $c_1 = (\tilde{\mu} K_p/4)^{2/3} - b - K_n$. Condition (10) implies that $c_1$ is positive, and (19), (20) and Lemma 3.4 applied to (12) imply that for $x \in [-L, -n] \cup [-l, l] \cup [n, L]$,

$$\frac{\tilde{\mu}}{2}\left(\frac{1}{\hat{\pi}_x} - (2L+1)\right) = \overline{m}_{\hat{\pi}} - m_x(\hat{\pi}) \geq c_1.$$

Therefore for $x \in [-L, -n] \cup [-l, l] \cup [n, L]$,

$$\hat{\pi}_x \leq \frac{1}{(2c_1/\tilde{\mu}) + 2L + 1} \leq \frac{\tilde{\mu}}{2((\tilde{\mu} K_p/4)^{2/3} - b - K_n)}$$

and the proof is complete. $\square$

3.2.2. *Small $\mu$.* Finally, we turn to the case of small $\mu$, where we only consider local maxima of $\overline{m}_\pi$ and then we can use the perturbation result Propositions 2.2 and 2.3 to say something about $V_\pi$.

We first establish Theorem 2.9, which says that if there is very little mass outside the interval $(-\lfloor M/2 \rfloor, \lfloor M/2 \rfloor)$, then the mass inside the interval $(-\lfloor M/2 \rfloor, \lfloor M/2 \rfloor)$ is concentrated at site 0.

PROOF OF THEOREM 2.9. Since $2q - 1 < M$, $B_{y,z} = b$ for $y, z \in (-q, q)$, therefore for $x \in (-q, 0) \cup (0, q)$,

$$
\begin{align*}
m_0(\hat{\pi}) &= K_0 \sum_z B_{-z} K_z \pi_z \geq b \sum_{z=-q+1}^{q-1} K_z \pi_z, \\
m_x(\hat{\pi}) &= K_x \left( b \sum_{z=x-M+1}^{x+M-1} K_z \pi_z + \sum_{z=-L}^{x-M} K_z \pi_z + \sum_{z=x+M}^{L} K_z \pi_z \right) \\
&\leq K_x \left( b \sum_{z=-q+1}^{q-1} K_z \pi_z + \sum_{z=-L}^{-q} K_z \pi_z + \sum_{z=q}^{L} K_z \pi_z \right) \\
&\leq K_1 \left( \varepsilon + b \sum_{z=-q+1}^{q-1} K_z \pi_z \right)
\end{align*}
$$

and

$$m_0(\hat{\pi}) - m_x(\hat{\pi}) \geq b(1-K_1) \sum_{z=-q+1}^{q-1} K_z \pi_z - K_1 \varepsilon$$

(21)
$$\geq b(1-K_1)K_{q-1}(1-\varepsilon) - K_1\varepsilon$$

$$= c > 0$$

if $\varepsilon < b(1-K_1)K_{q-1}/(b(1-K_1)K_{q-1} + K_1)$. Let $J = (-q, 0) \cup (0, q)$, then Lemma 2.4(b) says that

$$\frac{\sum_{x \in J} m_x(\hat{\pi})\hat{\pi}_x}{\sum_{x \in J} \hat{\pi}_x} + \frac{\tilde{\mu}}{2\sum_{x \in J} \hat{\pi}_x} = m_0(\hat{\pi}) + \frac{\tilde{\mu}}{2\hat{\pi}_0},$$

therefore

$$\frac{1}{\sum_{x \in J} \hat{\pi}_x} = \frac{2}{\tilde{\mu}} \left( \frac{\tilde{\mu}}{2\hat{\pi}_0} + m_0(\hat{\pi}) - \frac{\sum_{x \in J} m_x(\hat{\pi})\hat{\pi}_x}{\sum_{x \in J} \hat{\pi}_x} \right)$$

$$\geq \frac{1}{\hat{\pi}_0} + \frac{2}{\tilde{\mu}} \left( m_0(\hat{\pi}) - \max_{x \in J} m_x(\hat{\pi}) \right) \geq \frac{2c}{\tilde{\mu}}$$

by (21), which implies the desired conclusion. □

Now we focus on the local maxima that has their support spread $M$ sites apart, which seem to be the only local maxima from simulation. Suppose a subset $I$ of $E = [-L, L] \cap \mathbb{Z}$ has the properties that all $x \in I$ are at least $M$ apart. Define

$$\Delta^I = \{\pi \in \Delta : \pi_x = 0 \text{ for } x \notin I\}.$$

We observe that $\Delta^I$ is a closed subset of $\Delta$. We first establish Proposition 3.6 below that states a condition necessary for a local maximum in $\Delta^I$ to be a local maximum in $\Delta$, then establish Proposition 2.10, which has some results regarding the kinds of local maxima that various fitness functions can have.

LEMMA 3.5. *Let $k \in \mathbb{R}^+$ and*

$$k\Delta^I = \left\{ (\pi_{-L}, \ldots, \pi_L) : \pi_x \geq 0 \ \forall x \in I, \pi_x = 0 \ \forall x \notin I \text{ and } \sum_{x=-L}^{L} \pi_x = k \right\},$$

*then $\tilde{\pi} \in \Delta^I$ is the unique local maximum of $\overline{m}_\pi$ for $\pi$ lying in a small neighborhood in $\Delta^I$ if and only if $k\tilde{\pi} \in k\Delta^I$ is the unique local maximum of $\overline{m}_\pi$ for $\pi$ lying in a small neighborhood in $k\Delta^I$.*

PROOF. If $\tilde{\pi} \in \Delta^I$ is the unique local maximum of $\overline{m}_\pi$ for $\pi$ lying in a small neighborhood in $\Delta^I$, then the Hessian matrix of $\overline{m}_\pi$ at $\tilde{\pi}$ is positive definite, and Lemma 3.1 implies that $m_x(\hat{\pi})$ are all equal for $x \in I$. Thus $m_x(k\tilde{\pi}) = K_x \sum_{z \in I} B_{x-z} K_z k\tilde{\pi}_z = k m_x(\tilde{\pi})$ are all equal for $x \in I$ as well. This shows $k\tilde{\pi}$ is a local extremum of $\overline{m}_\pi$ for $\pi \in k\Delta^I$. To verify it is a local maximum, we define $I' = I \setminus \{p\}$ where $p$ is an arbitrary member of $I$, rewrite $\overline{m}_\pi$ in terms of $x \in I'$, and calculate its first and second derivatives:

$$\overline{m}_\pi = \sum_{x,z \in I'} K_x \pi_x B_{x-z} K_z \pi_z + 2 K_p \left(k - \sum_{x \in I'} \pi_x\right) \sum_{z \in I'} B_{p-z} K_z \pi_z$$
$$+ B_0 K_p^2 \left(k - \sum_{x \in I'} \pi_x\right)^2,$$

$$\frac{\partial \overline{m}_\pi}{\partial \pi_w} = 2 K_w \sum_{z \in I'} B_{w-z} K_z \pi_z - 2 K_p \sum_{z \in I'} B_{p-z} K_z \pi_z$$
$$+ 2 K_p \left(k - \sum_{x \in I'} \pi_x\right) B_{p-w} K_w - 2 B_0 K_p^2 \left(k - \sum_{x \in I'} \pi_x\right),$$

$$\frac{\partial^2 \overline{m}_\pi}{\partial \pi_w \partial \pi_y} = 2 K_w B_{w-y} K_y - 2 K_p B_{p-y} K_y - 2 K_p B_{p-w} K_w + 2 B_0 K_p^2,$$

where $w, y \in I'$. We observe that the second derivatives do not depend on $k$, therefore the Hessian matrix of $\overline{m}_\pi$ is also positive definite at $k\tilde{\pi}$, and $k\tilde{\pi}$ is the unique local maximum lying in a small neighborhood in $k\Delta^I$. The proof of the reverse direction is similar. □

PROPOSITION 3.6. *Suppose $\tilde{\pi} \in \Delta^I$ is the unique local maximum of $\overline{m}_\pi$ for $\pi$ lying in a small neighborhood of $\tilde{\pi}$ in $\Delta^I$, $\tilde{\pi}_x > 0$ for all $x \in I$, and $m_x(\tilde{\pi}) = m_1$ for all $x \in I$.*

*(a) If $m_x(\tilde{\pi}) \leq m_2 < m_1$ for all $x \notin I$, then $\tilde{\pi}$ is also the unique local maximum of $\overline{m}_\pi$ for $\pi$ lying in a sufficiently small neighborhood in $\Delta$.*

*(b) If there exists $y \in E \setminus I$ where $m_y(\tilde{\pi}) = m_3 > m_1$, then $\tilde{\pi}$ is not a local maximum of $\overline{m}_\pi$ for $\pi \in \Delta$.*

REMARK 3.7. If the set $I$ consists of a singleton $y$, then $\delta_y \in \Delta^I$ is trivially the unique local maximum of $\overline{m}_\pi$ for $\pi$ lying in a small neighborhood of $\tilde{\pi}$ in $\Delta^I$, which is an empty set. In this case, to verify that $\delta_y$ is also the unique local maximum of $\overline{m}_\pi$ for $\pi$ lying in a sufficiently small neighborhood in $\Delta$, we only need to check that $m_x(\delta_y) < m_y(\delta_y)$ for all $x \neq y$.

PROOF OF PROPOSITION 3.6. (a) We examine the fitness of $\pi$ in a small neighborhood

$$A = \left\{\pi \in \Delta : \max_{x \in E} |\pi_x - \tilde{\pi}_x| < \varepsilon\right\}$$

of $\tilde{\pi}$ in $\Delta$:

$$\begin{aligned}
\overline{m}_\pi &= \sum_{x \in I, z \in I} K_x \pi_x B_{x-z} K_z \pi_z + \sum_{x \in I, z \notin I} K_x \pi_x B_{x-z} K_z \pi_z \\
&\quad + \sum_{x \notin I, z \in E} K_x \pi_x B_{x-z} K_z \pi_z \\
&\leq \sum_{x \in I, z \in I} K_x \pi_x B_{x-z} K_z \pi_z + 2 \sum_{x \notin I, z \in E} K_x \pi_x B_{x-z} K_z \pi_z \\
&= \sum_{x \in I, z \in I} K_x \pi_x B_{x-z} K_z \pi_z + 2 \sum_{x \notin I} \pi_x m_x(\pi).
\end{aligned} \quad (22)$$

Let $c = \sum_{x \notin I} \pi_x$. By Lemma 3.5, the configuration $\tilde{\pi}_1 \in (1-c)\Delta^I$ that maximizes (locally) the first term on the right-hand side of (22) is $(1-c)\tilde{\pi}$, which means that it satisfies

$$\sum_{x \in I, z \in I} K_x \pi_x B_{x-z} K_z \pi_z \leq (1-c)^2 m_1$$

if $\varepsilon$ in the definition of $A$ is sufficiently small. If $\varepsilon$ is sufficiently small, then because $m_x(\pi)$ is a continuous function of $\pi$ for all $x$, we have

$$m_x(\pi) \leq m_2 + \frac{m_1 - m_2}{2}$$

for $x \notin I$. Applying the two estimates above to (22), we obtain

$$\overline{m}_\pi \leq (1-c)^2 m_1 + 2c\left(m_2 + \frac{m_1 - m_2}{2}\right) = m_1 - c(m_1 - m_2) + c^2 m_1 < m_1$$

if $c \leq (2L+1)\varepsilon$ is strictly positive but sufficiently small. Hence $\tilde{\pi}$ is also the unique local maximum of $\overline{m}_\pi$ for $\pi$ lying in a small neighborhood in $\Delta$.

(b) Let $w \in I$, then $\tilde{\pi}_w > 0$ and $m_w(\tilde{\pi}) = m_1$. Along the line $\tilde{\pi} + p(\delta_y - \delta_w)$,

$$\begin{aligned}
\left.\frac{\partial \overline{m}_{\pi + p(\delta_y - \delta_w)}}{\partial p}\right|_{\pi = \tilde{\pi}, p = 0} &= \left.\frac{\partial}{\partial(\pi_y + p)}\right|_{\pi = \tilde{\pi}, p = 0} \overline{m}_{\pi_x + p(\delta_{y,x} - \delta_{w,x})} \frac{\partial(\pi_y + p)}{\partial p} \\
&\quad + \left.\frac{\partial}{\partial(\pi_w - p)}\right|_{\pi = \tilde{\pi}, p = 0} \overline{m}_{\pi_x + p(\delta_{y,x} - \delta_{w,x})} \frac{\partial(\pi_w - p)}{\partial p} \\
&= 2(m_y(\tilde{\pi}) - m_w(\tilde{\pi})) \\
&= 2(m_3 - m_1),
\end{aligned}$$

which is strictly positive. Therefore $\tilde{\pi}$ is not a local maximum of $\overline{m}_\pi$ for $\pi \in \Delta$. □

Now we use the above result to establish Proposition 2.10, which explicitly computes some local maxima, when the dimension of $\Delta^I$ is low enough (less than 3) to enable us to do hands-on computation.

PROOF OF PROPOSITION 2.10. (a) By Lemma 3.1, at any local maximum $\hat{\pi}$, the support of $\hat{\pi}$ must have equal fitness. We will show that either $x_j + 1$ or $x_{j+1} - 1$ is more fit than sites in $\{x_1, \ldots, x_k\}$. If $x_{j+1} \geq 1$, then since $K$ is unimodal by assumption, $K_{x_{j+1}-1} > K_{x_{j+1}}$. Otherwise, $x_{j+1} \leq 0$ and $x_j \leq -1$, therefore $K_{x_j+1} > K_{x_j}$. So either $K_{x_{j+1}-1} > K_{x_{j+1}}$ or $K_{x_j+1} > K_{x_j}$ or both. If $K_{x_j+1} > K_{x_j}$ (the case of $K_{x_{j+1}-1} > K_{x_{j+1}}$ is similar), then

$$m_{x_j+1}(\hat{\pi}) = K_{x_j+1}\left(\sum_{z \in \{x_1,\ldots,x_j\}} B_{x_j+1-z}K_z\hat{\pi}_z + \sum_{z \in \{x_{j+1},\ldots,x_k\}} B_{x_j+1-z}K_z\hat{\pi}_z\right).$$

Since $x_{j+1} - (x_j + 1) \geq M$, all $B_{x_j+1-z}$ in the second sum above are 1 and equal to $B_{x_j-z}$ for $z \in \{x_{j+1},\ldots,x_k\}$ and $B_{x_j+1-z} \geq B_{x_j-z}$ for $z \in \{x_1,\ldots,x_j\}$ in the first sum above, therefore

$$m_{x_j+1}(\hat{\pi}) \geq K_{x_j+1}\left(\sum_{z \in \{x_1,\ldots,x_j\}} B_{x_j-z}K_z\hat{\pi}_z + \sum_{z \in \{x_{j+1},\ldots,x_k\}} B_{x_j-z}K_z\hat{\pi}_z\right)$$

$$> K_{x_j}\left(\sum_{z \in \{x_1,\ldots,x_j\}} B_{x_j-z}K_z\hat{\pi}_z + \sum_{z \in \{x_{j+1},\ldots,x_k\}} B_{x_j-z}K_z\hat{\pi}_z\right)$$

$$= m_{x_j}(\hat{\pi}),$$

therefore $\hat{\pi}$ cannot be a local maximum of $\overline{m}_\pi$ by Proposition 3.6(b).

(b) The fact that $\delta_0$ is a stationary point of (9) is obvious; in fact, it holds for any $K$ and $b$. But we need to check it is a local maximum of $\overline{m}_\pi$ if $K_M < b$. We compute the fitness of all sites in $E$ when $\pi = \delta_0$. For $x \in [-M+1, M-1]$,

$$m_x(\pi) = \begin{cases} bK_x, & \text{if } x \in [-M+1, M-1], \\ K_x, & \text{if } x \in [-L, -M] \cup [M, L]. \end{cases}$$

If $K_M < b$, then $m_x < m_0$ for all $x \neq 0$ since $K$ is increasing in $[-L, 0]$ and decreasing in $[0, L]$. Proposition 3.6(a) and Remark 3.7 imply that $\delta_0$ is a local maximum of $\overline{m}_\pi$.

Now we deal with the case of $K_M > b$. Since $m_M(\delta_0) = K_M > b = m_0(\delta_0)$, Proposition 3.6(b) and Remark 3.7 imply that $\delta_0$ is not a local maximum of $\overline{m}_\pi$.

(c) For $\hat{\pi} = p\delta_{-x} + (1-p)\delta_{-x+M}$ where $p$ is defined in the statement, brute force calculation shows that

$$m_{-x}(\hat{\pi}) = m_{-x+M}(\hat{\pi}) = \frac{K_{-x}^2 K_{-x+M}^2 (1-b^2)}{2K_{-x}K_{-x+M} - bK_{-x}^2 - bK_{-x+M}^2}$$

and

$$\frac{\partial^2 \overline{m}_{\hat{\pi}}}{\partial p^2} = bK_{-x}^2 + bK_{-x+M}^2 - 2K_{-x}K_{-x+M},$$

which is $< 0$ if $bK_{-x}^2 + bK_{-x+M}^2 < 2K_{-x}K_{-x+M}$. This verifies that $\hat{\pi}$ is the unique local maximum of $\overline{m}_\pi$ for $\pi$ lying in a small neighborhood of $\hat{\pi}$ in $\Delta^{\{-x,-x+M\}}$.

It remains to check that all sites other than $-x$ and $-x+M$ are less fit, which calculations of fitness at these sites show to be true if $b$, $K_{-x-M}$, and $K_{-x+2M}$ are all $< K_{-x}K_{-x+M}(1+b)/(K_{-x} + K_{-x+M})$.

(d) This result can be proved by brute force calculation, just like part (c). We omit the details. $\square$

**Acknowledgments.** The author is grateful to Martin Barlow, Michael Doebeli, Edwin Perkins and John Walsh for valuable discussions during this research.

Department of Statistics
University of Oxford
1 South Parks Road
Oxford OX1 3TG
United Kingdom
E-mail: fyz@stats.ox.ac.uk
URL: http://www.stats.ox.ac.uk/~fyz